\documentclass[11pt]{amsart}
\usepackage{amsmath,amssymb,latexsym,soul,cite,amsthm,color,enumitem,graphicx,tikz, mathtools,microtype}

\usepackage{amsfonts,amsmath,amsthm,latexsym,amssymb,marvosym,esint,xfrac,bbm,amssymb}

\usepackage{cite}

\usepackage{graphics, epsfig}
\tikzstyle{mybox} = [draw=black, very thick, rectangle, rounded corners, inner ysep=5pt, inner xsep=5pt]
\usepackage{appendix}

\usepackage[colorlinks=true, urlcolor=blue, citecolor=red, linkcolor=blue]{hyperref}

\usepackage[english]{babel}
\usepackage[left=2.5cm,right=2.5cm,top=2.7cm,bottom=2.7cm]{geometry}

\numberwithin{equation}{section}

\DeclareMathOperator{\Div}{div}
\newtheorem{theorem}{Theorem}[section]
\theoremstyle{plain}
\newtheorem{lemma}[theorem]{Lemma}
\theoremstyle{plain}

\theoremstyle{plain}

\newtheorem{definition}[theorem]{Definition}
\theoremstyle{definition}
\newtheorem{remark}[theorem]{Remark}

\newcommand{\N}{{\mathbb N}}

\newcommand{\R}{{\mathbb R}}

\newcommand{\e}{\epsilon}

\newcommand{\beq}{\begin{equation}}
\newcommand{\eeq}{\end{equation}}
\renewcommand{\le}{\leqslant}
\renewcommand{\ge}{\geqslant}

\newcommand{\C}{\mathcal{C}}

\def\Xint#1{\mathchoice
    {\XXint\displaystyle\textstyle{#1}}%
    {\XXint\textstyle\scriptstyle{#1}}%
    {\XXint\scriptstyle\scriptscriptstyle{#1}}%
    {\XXint\scriptscriptstyle\scriptscriptstyle{#1}}%
    \!\int}
\def\XXint#1#2#3{\setbox0=\hbox{$#1{#2#3}{\int}$}
    \vcenter{\hbox{$#2#3$}}\kern-0.5\wd0}
\def\bint{\Xint-}
\def\bint{\Xint-}
\def\dashint{\Xint{\raise4pt\hbox to7pt{\hrulefill}}}
\def\dashiint{\bint\kern-0.15cm\bint}

\def\tr{(u-k)_{-}}
\def\A{\mathbb{A}}
\def\Q{\mathcal{Q}}
\def\K{\mathcal{K}}
\def\dive{\mathrm{div}}
\def\d{\mathrm{d}}
 \def\B{\mathcal{B}}
\def\dim{\mathrm{dim}}
\def\dist{\mathrm{dist}}
\def\H{\mathcal{H}}
\def\I{\mathcal{I}}
\def\L{\mathcal{L}}
\def\C{\mathcal{C}}
\def\P{\mathcal{P}}
\def\T{\mathcal{T}}
\def\S{\mathcal{S}}

\title[Integral Harnack Estimates and Extinction Profile] {The impact of intrinsic scaling on the rate of extinction for anisotropic non-newtonian fast diffusion}
\author[Ciani, Henriques, Skrypnik]{Simone Ciani {\&} Eurica Henriques {\&} Igor I. Skrypnik}
\address{Department of Mathematics of the University of Bologna, Piazza Porta San Donato, 5, 40126 Bologna, Italy}
\email{simone.ciani3@unibo.it}
\address{Centro de Matem\'atica, Universidade do Minho - Polo CMAT-UTAD Departamento de Matem\'atica Universidade de Tr\'as-os-Montes e Alto Douro, Vila Real, Portugal}
\email{eurica@utad.pt}
\address{Institute of Applied Mathematics and Mechanics, National Academy of Sciences of Ukraine, Gen. Batiouk Str. 19, 84116 Sloviansk, Ukraine}
\email{ihor.skrypnik@gmail.com}

\begin{document}
\begin{abstract} 
We study the decay towards the extinction that pertains to local weak solutions to  fully anisotropic equations whose prototype is  
\begin{equation*} 
    \partial_t u= \sum_{i=1}^N \partial_i (|\partial_i u|^{p_i-2} \partial_i u), \qquad 1<p_i<2.
\end{equation*} \noindent Their rates of extinction are evaluated by means of several integral Harnack-type inequalities which constitute the core of our analysis and that are obtained for anisotropic operators having full quasilinear structure. Different decays are obtained when considering different space geometries. The approach is motivated by the research of new methods for strongly nonlinear operators, hence dispensing with comparison principles, while exploiting an intrinsic geometry that affects all the variables of the solution. \vskip0.2cm \noindent 

\end{abstract}

\maketitle
	\begin{center}
		\begin{minipage}{9cm}
			\small
   \tableofcontents
		\end{minipage}
	\end{center}



\def\tr{(u-k)_{-}}
\def\R{\mathbb{R}}
\def\A{\mathbb{A}}
\def\Q{\mathcal{Q}}
\def\b{\mathbb{B}}
\def\K{\mathcal{K}}
\def\k{\mathbb{K}}
\def\E{\mathcal{E}}
\def\q{\mathbb{Q}}

\def\N{\mathbb{N}}
\def\dive{\mathrm{div}}
\def\d{\mathrm{d}}
 \def\B{\mathcal{B}}
\def\dim{\mathrm{dim}}
\def\dist{\mathrm{dist}}
\def\H{\mathcal{H}}
\def\I{\mathcal{I}}
\def\L{\mathcal{L}}
\def\C{\mathcal{C}}
\def\T{\mathcal{T}}
\def\S{\mathcal{S}}
\def\P{\mathcal{P}}
\def\e{\epsilon}


\addtocontents{toc}{\protect\setcounter{tocdepth}{1}}

\section{Introduction}
\noindent For an open bounded set $\Omega \subset \R^N$ and a positive time $T$, we consider anisotropic differential equations whose prototype is the following
\begin{equation}\label{Prototype}
\partial_t u-\Delta_{{\bf p}} u :=\partial_t u-\sum_i \partial_i ( |\partial_i u |^{p_i-2} \partial_i u )=0, \qquad \text{weakly in} \qquad \Omega_T= \Omega \times [0,T].
\end{equation}
\noindent Differential operators as $(\partial_t  - \Delta_{\bf p})$ above appear already in the seminal work \cite{Lions}, in the guise of the prototype example of operators obtained as the sum of monotone ones. They enjoy many interesting properties (see for instance the book \cite{AS-book2}) whose interpretation has led to a rich mathematical theory (see for instance \cite{IlinBesovNikolski}, \cite{Cianchi}, \cite{Tartar}, \cite{Troisi}). Nonetheless, even after more than half a century, the basic regularity properties of local weak solutions to equations \eqref{Prototype} remain an open problem (see for instance \cite{AMS}, \cite{Brasco}, \cite{CMV}). Besides the theoretical intrinsic interest and challenge, this kind of equations appear in various physical contexts (see Chap. IV of \cite{AS-book1}), unveiling the mathematical description of diffusion processes for which the propagation has a different non-Newtonian behavior along each coordinate axis; as well as  modeling electro-rheological fluids (see for instance the seminal paper \cite{Rajagopal} or the book \cite{Ruzika}), in particular when the stress tensor is a function of an electromagnetic field that varies on each coordinate direction. 
\vskip0.1cm \noindent 
This work is developed for the so-called  {\it fast diffusion} regime, $1 <p_i<2$ for all $i \in \{1,\dots,N\}$, which seems to unfold very strong properties of solutions. The precise attribute we are interested in is the property of {\it extinction in finite time} of local weak solutions to \eqref{Prototype}, meaning that  there exists a finite time $T^*<T$, called time of extinction, such that the solution $u$ vanishes out from $T^*$: 
\[ \exists \  T^* \in [0,T]: \qquad u(\cdot, t) \equiv 0 \ ,  \qquad  \forall t \ge T^*. \]

\noindent This property is enjoyed by the solutions to the parabolic $p$-Laplacean equation
\begin{equation}\label{plap}
\partial_t u -\Delta_p u := \partial_t u-\Div ( |\nabla u |^{p-2} \nabla u )=0, \qquad \text{weakly in} \qquad \Omega_T=\Omega \times [0,T],
\end{equation} \noindent and it affects preponderantly the nature and behavior of solutions (see \cite{DB} or, more in general, \cite{BenCra} and \cite{Diaz-anulacion}).\newline For instance, in \cite{DBKW} the authors  show that a point-wise Harnack inequality cannot be found for the solutions to \eqref{plap} in the sub-critical range $1 <p<2N/(N+1)$; while in the super-critical range $2N/(N+1)<p<2$ the phenomenon of expansion of positivity is closely related to the singular character of the operator, that privileges the elliptic behavior to the diffusive one, as soon as the modulus of ellipticity $|\nabla u |^{p-2} \nabla u$ blows up.

\noindent To the very interesting properties of singular equations, the operator \eqref{Prototype} adds the fascinating ones of anisotropy. In \cite{Vaz2}, the asymptotic behavior is studied through the analysis of self-similarity, showing that new mathematical methods need to be developed in order to overcome the strong non-uniqueness phenomena and to construct suitable barriers. In \cite{AS2}, the authors show that these anisotropic equations are, in a certain sense, richer than their $p$-Laplacean counterpart; indeed, for solutions to equations as \eqref{Prototype} within the more relaxed condition $1<p<2$ (here $p$ is an average of $p_i$s, see Section \ref{Preliminaries}) the dichotomy finite speed of propagation/extinction in finite time is no longer valid and it is replaced by conditions on the growth exponents $p_i$s taking into account the competition between diffusions.

\vskip0.2cm \noindent Solutions to singular $p$-Laplacean equations as \eqref{plap}, have a decay toward extinction (see \cite{DBKW}) that follows the law
\[ \|u(\cdot, t)\|_{\infty, B_{\rho}}  \leq \gamma \bigg( \frac{T^{*}-t}{\rho^p}\bigg)^{\frac{1}{2-p}}, \qquad \forall \rho,t>0\, : \, \,  B_{\rho}\times ((t+T^{*})/2, \, T^{*}] \subset \Omega_T,
\] being $B_{\rho}$ the ball of radius $\rho$ and $\gamma$ a positive constant depending only on the data $\{N,p\}$. In the present work we show that the decay profile of extinction of solutions to equations of the kind of \eqref{Prototype} is the same as the one to the $p$-Laplacean if one considers a particular space-geometry, 
\[\|u(\cdot, t)\|_{\infty, \K_{\rho}(T^*-t)}\leq \gamma \bigg( \frac{T^*-t}{\rho^p}\bigg)^\frac{1}{2-p}, \qquad \forall \rho,t>0 \, : \, \, \K_{\rho}(T^*-t) \times ((t+T^{*})/2, \, T^{*}] \subset \Omega_T,
\] being $\gamma$ a positive constant depending only on the data, and, for any fixed $\tau>0$ 
\begin{equation} \label{krhot} \K_{\rho}(\tau)= \prod_i \bigg{\{}|x_i|< \rho^{\frac{p}{p_i}} \bigg( \frac{\tau}{\rho^p} \bigg)^{\frac{p-p_i}{p_i(2-p)}}  \bigg{\}}, \qquad \text{being} \qquad p=N/\bigg(\sum_i 1/p_i\bigg).\end{equation}

\noindent This particular space geometry, which we refer to as {\it intrinsic geometry} (see Section \ref{main-results}), has interesting features: although the cylinder $\K_{\rho}(T^*-t)$ degenerates in these directions $x_i$ for which $p_i>p$ when $t$ approaches $T^*$, it preserves its volume regardless of the time level undertaken; and more, when $p_i \equiv p$ for all $i =1, \dots, N$, the set $\K_{\rho}(\tau)$ is the classical cube.

\noindent We also show that the decay rate of a solution $u$ to equations of the type \eqref{Prototype} can be estimated within a geometry that is non-degenerative, but at the price of a more complex rate
\[\|u(\cdot, t) \|_{\infty, \k_{\rho}} \leq \gamma \sum_i \bigg(\frac{T^*-t}{\rho^p} \bigg)^{\frac{\lambda_i}{(2-p_i)\lambda}}, \qquad \text{being} \quad \lambda_i=N(p_i-2)+p,\] $\lambda=N(p-2)+p$ (as usual) and $\gamma$ a positive constant depending on the data. Here the geometry will be referred to as the {\it standard geometry}, being based on cubes as
\begin{equation}\label{anisociccio}
    \k_{\rho}= \prod_i \bigg{\{} |x_i|< \rho^{\frac{p}{p_i}}\bigg{\}}, \qquad \qquad \rho>0.
\end{equation} Unlike the {\it intrinsic geometry} considered before, this one does not take into account the time variable. Again, when $p_i \equiv p$ for all $i=1,\dots, N$, the set $\k_{\rho}$ is the classic cube of hedge $2 \rho$. It is clear that the extinction rate in this case will depend on the smallness of $T^*-t$ and the maximum of the exponents in the sum.

\vskip0.2cm \noindent
It is the precise aim of our study to carry out an analysis of these two rates of extinction within these two different underlying geometries. The method of derivation of these decay rates has its own mathematical interest: confirming the well-known principle that the run itself can be more instructive than the final destination, we obtain the above behaviour of solutions from various Harnack-type estimates.  These inequalities are found in three different topologic settings: $L^1_{loc}(\Omega)$, $L^1_{loc}(\Omega)$-$L^{\infty}_{loc}(\Omega)$ and $L^r_{loc}(\Omega)$-$L^r_{loc}(\Omega)$ backward in time,  and all of them are new for solutions to operators as \eqref{Prototype} (we refer to Section \ref{main-results} for the precise statements).\newline \noindent Here below we give an example of what we mean by Harnack-type estimates in the $L^1_{loc}(\Omega)$-topology, or, in short, $L^1$-$L^1$ Harnack-type inequality.

\vskip0.2cm 

\noindent \hspace*{0mm}\begin{tikzpicture}
\node [mybox] (box){%
    \begin{minipage}{.96\textwidth}
\small{{\it {\bf{$L^1$-$L^1$ Harnack-type inequality}}}}\vskip0.2cm \noindent 
Let $u$ be a non-negative local weak solution to \eqref{Prototype} in $\R^N \times \R_0^+$ and let $\rho,t$ be positive fixed numbers.\newline
Then, the following two estimates hold true in their respective space configurations.
\vskip0.3cm \noindent 
\begin{itemize}
    \item[1] Let $\K_{\rho}(t)$ be defined as 
    in \eqref{krhot}. Then there exists a constant $\gamma(N,p_i)>1$ such that 
    \begin{equation*} \label{francesco}
    \sup_{0\leq \tau \leq t} \int_{\K_{\rho}(t)} u(x,\tau)\, dx \leq \gamma \inf_{0\leq\tau\leq t} \int_{2\K_{\rho}(t)} u(x,\tau) \, dx + \gamma \bigg( \frac{t}{\rho^\lambda} \bigg)^{\frac{1}{2-p}}.
\end{equation*}

\vskip0.4cm \noindent
    \item[2]  Let $\k_{\rho}$ be defined as in \eqref{anisociccio}. Then there exists a constant $\gamma(N,p_i)>1$ such that  \[
     \sup_{0\leq \tau \leq t} \int_{\k_{\rho}} u(x,\tau)\, dx \leq \gamma \inf_{0\leq \tau \leq t} \int_{2\k_{\rho}} u(x,\tau) \, dx + \sum_{i} \left(\frac{t}{\rho^{\lambda_i}}\right)^{\frac{1}{2-p_i}} .\] 
\end{itemize}
    \end{minipage}
};
\end{tikzpicture}%
\newline \noindent

\subsection*{Novelty and Significance}
\subsubsection*{Origins} To the best of our knowledge, the idea of a Harnack-type estimate in the topology of $L^1_{loc}(\Omega)$ had its first appearance in \cite{DBH-sing} for the prototype $p$-Laplacean equation, and it was used in \cite{DBKW} with the aim of giving a bound from below to its solutions in a small cylinder, so to prove a point-wise Harnack inequality. There these integral Harnack-type estimates are first used to evaluate the time of extinction of solutions.\newline
The method has been reported in (\cite{DB}, Chap. VII) for solutions to the prototype singular equation ($1<p<2$). A proof for $p$-Laplacean type equations with full quasilinear structure can be found first in the paper \cite{Annali} and then in the monograph \cite{DBGV-mono}, again with the aim of obtaining a bound from below toward the determination of a point-wise Harnack-type inequality.\vskip0.1cm \noindent 
All these estimates are unknown for anisotropic equations such as \eqref{Prototype}. In contrast with the few results available in literature (see for instance \cite{CMV}, \cite{Vaz2}) that use crucially the invariance and comparison properties of the prototype equation, we derive here the aforementioned Harnack-type inequalities for the full-quasilinear structure operator (see definition \eqref{EQ-weak}-\eqref{structure-conditions}) adopting a technique that dispenses with comparison principles and treats equations that have bounded and measurable coefficients. For this whole spectrum of equations we derive the decay rate of extinction.\vskip0.2cm \noindent 
As anticipated, in the {\it cours d'oevre} for the evaluation of the extinction rate, we derive backward $L^r_{loc}(\Omega)$-$L^{\infty}_{loc}(\Omega)$ estimates that have their own mathematical interest (see Theorems \ref{lr-decay}, \ref{elr-decay}). For their derivation, we assume that the solutions are locally bounded: this is a crucial point for the regularity theory of anisotropic $p$-Laplacean equations, as a condition on the spareness of the exponents $p_i$s is necessary already for the elliptic case (see for instance \cite{Giaquinta},\cite{Marcellini}). From the (anisotropic) parabolic point of view, the theory of local boundedness is reasonably complete, see for instance \cite{CVV}, \cite{DMV}, \cite{LY}. Finally, these $L^r_{loc}(\Omega)$-$L^{\infty}_{loc}(\Omega)$ estimates are reminiscent of the isotropic case (see for instance \cite{DBKW}) and are obtained through the successive application of standard $L^r_{loc}(\Omega)$-$L^{\infty}_{loc}(\Omega)$ estimates (Theorems \ref{lrlinfestimate}, \ref{sup}) with backwards $L^r_{loc}(\Omega)$ ones (see Theorems \ref{Y}, \ref{lrestimate}). We refer to \cite{DBGV-mono} and the references therein for the isotropic counterpart.
\vskip0.1cm 

\noindent The lack of (known) regularity of solutions encumbers the research for applications on models directly intertwined with \eqref{Prototype} (see \cite{AS-book1} Chap. IV). Nonetheless, these operators reveal a very interesting picture of the underlying nonlinear analysis and competitive behaviour between different diffusions. 

    \subsubsection*{The role of intrinsic geometry} A satisfying study of anisotropic operators as \eqref{Prototype} cannot be brought on regardless of the self-similar geometry embodied in the operator itself. This is already understood in the case of the evolutionary $p$-Laplacean equation, where has been shown that a Harnack inequality holds true only in a particular geometry, called {\it intrinsic geometry}. We refer to \cite{DB} and \cite{Urbano} for insights on this topic. Roughly speaking, in the regularity theory of diffusive $p$-Laplacean equations, time is linked to space by a relation that takes into account the solution itself, as $t= \rho^p u_o^{2-p}$, supposing $u_o>0$ is the value of the solution at a point. In the case of anisotropic operators behaving like \eqref{Prototype}, the full power of self-similar geometry is needed, and the scaling factor depending on $u_o$ enters also the in space variables. As a concrete example, in the degenerate case and for solutions $u$ of \eqref{Prototype} in $S_{\infty}=\R^N \times \R_+$, a point-wise Harnack inequality takes the following form (we refer to \cite{CMV}):

\begin{equation*}\label{isoHarnack}
  \frac{1}{\gamma}\sup_{\K_{\rho}(M)}u(\,\cdot\, , - M^{2- p}\, (C_{2}\, \rho)^{p} )\le  u_0 \le \gamma \inf_{\K_{\rho}(M)} u(\,\cdot\, ,   M^{2-p}\, (C_{2}\, \rho)^{{p}})
    \end{equation*}
    with $M= (u_o/C_1)$, being $\gamma, C_1,C_2 $ positive constants depending only on $\{N,p_i\}$ . In the available literature, $L^1$-$L^1$ Harnack-type estimates are derived for the diffusive $p$-Laplacean operators (see \cite{DBGV-mono}) without the use of a particular intrinsic geometry. Here we overcome the difficulty of the non-homogeneity of the operator by setting an intrinsic geometry that depends also on time, as $\K_{\rho}(t)$ in \eqref{krhot}, which considers self-similar space-cubes as
    \[\K_{\rho}(M)= \prod_i \bigg{\{}  |x_i|< \rho^{\frac{p}{p_i}}M^{\frac{p_i-p}{p_i}}\bigg{\}}, \quad \text{with} \quad M= \bigg( \frac{t}{\rho^p}\bigg)^{\frac{1}{2-p}}.\]
    \noindent In this case, the particular self-similar factor $M$ depends on the radius and on the {\it a priori} chosen time level $t$, and has the interesting feature of reestablishing the homogeneity in the estimates. With a little abuse of notation, along the text we still call this geometry {\it intrinsic geometry}, because the quantity $M$ here above is always related to some norm of $u$ in applications (see for instance the use of \eqref{kappa} and \eqref{nu}). \vskip0.1cm \noindent 
    A last word in honor of the standard geometry $\k_{\rho}$ is due. Local integral $L^1$-$L^{\infty}$ Harnack-type inequalities hold true also in this case (see Theorems \ref{el1-linfty-THM}-\ref{el1l1}), which is when one considers $M=1$; but the anisotropy is inevitably carried over into a sum of the quantities $t/\rho^{p}$ on the right-hand side of the estimates, with different powers depending on $p_i$s. A novel method is also used in this case, which we believe to be useful also for other nonlinear operators.

    \subsubsection*{Applications and Future Perspectives}
    The range of application of the Harnack-type inequalities we are about to describe is very wide. As for the main purpose of the present work, they can be used to estimate the decay of the solution at the extinction time; and, assuming an integrable initial datum $\|u_0\|_{L^1(\R^N)}$ they imply a certain conservation of the mass of the solution in time. \newline In addition, not only these Harnack-type estimates are very important for the convergence of approximating solutions when dealing with the problem of the existence (see for instance \cite{DBH-sing}), but also they proved to be useful to control the measure of level sets and to give a short proof of solutions' H\"older continuity (see for instance \cite{CV-Rendiconti} for the isotropic case).

\subsection*{Method} The  Harnack-type estimates that are obtained throughout the paper, for each one of the mentioned geometries, have as common starting point some general energy estimates, that are collected in the Appendix. Although these energy estimates are non-trivial, they are similar to the isotropic ones (see Section \ref{appendix}); hence we decided to postpone their presentation so as to leave space to what is really new in the anisotropic context. 

\vskip0.2cm \noindent 
Our first step is to derive $L^1$-$L^1$ Harnack-type estimates by means of testing the equation with negative powers of the solution and a combined nonlinear iteration. In a second step, we study the $L^r$-$L^\infty$ inequalities by suitably adapting the classic De Giorgi-Moser scheme; here we use the $L^r$-norm of the solution chained with the energy estimates provided by the equation in a certain geometry. Finally, we nest these inequalities with a backward $L^r$ estimate to derive $L^r$-$L^{\infty}$ inequalities in terms of the initial datum $u_0$; combining these with the first obtained $L^1$-$L^1$ estimates we derive the $L^{1}$-$L^{\infty}$ Harnack-type estimates given by  Theorems \ref{l1-linfty-THM}, \ref{el1-linfty-THM}.

\subsection*{Structure of the paper} In Section \ref{main-results}, we define the anisotropic operators with full quasilinear structure and state the main Theorems. Then, in Section \ref{Preliminaries}, we give the definition of local weak solution and the proper functional spaces for it; along with the main notation used throughout the paper. In Section \ref{l1l1proofs}, we present the proofs of the first two Theorems, both concerning $L^1$-$L^1$ Harnack-type estimates, but specializing the geometry in each case. In a similar fashion, in Section \ref{lrlinfty}, we provide the proofs of the backward $L^r$-$L^{\infty}$ estimates, again distinguishing the two geometries. Finally, short Section \ref{section-l1linfty} concludes with the main Theorems, while the last Section, Appendix \ref{appendix}, presents the main energy estimates used along our analysis and some standard iteration Lemmata.

\section{Main Results and Applications} \label{main-results}
\noindent We consider singular parabolic nonlinear partial differential equations of the form
\begin{equation}\label{EQ}
    \partial_t u - \dive{A(x,t,u,Du)}=B(x,t,u,Du), \qquad \text{weakly in}\, \, \Omega_T= \Omega \times [0,T],
\end{equation} \noindent where the functions $A=(A_1,\dots, A_N): \Omega_T \times \R^{N+1} \rightarrow \R^N$ and $B: \Omega_T \times \R^{N+1} \rightarrow \R$ are Caratheodory functions that satisfy the structure conditions, for $ 1<p_i<2$, for all $ i=1, \dots, N,$
\begin{equation}\label{structure-conditions}
\begin{cases}
A_i(x,t,s,\xi)\, \xi_i \ge C_o |\xi_i|^{p_i}-C^{p_i},\\
|A_i(x,t,s,\xi)|\leq C_1 |\xi_i|^{p_i-1}+C^{p_i-1},  \\
|B(x,t,s,\xi)|\leq  \displaystyle{\sum_i C\bigg(|\xi_i|^{p_i-1}+C^{p_i-1}\bigg)}\, , 
\end{cases}
\end{equation}\noindent for almost every $(x,t) \in \Omega_T$ and for all $(s,\xi) \in \R\times \R^N$, where $C_o, C_1$ are positive constants and $C$ is a non-negative constant that distinguishes between the cases when the equation to be homogeneous (when $C=0$) from when it is not.\newline We will say that a positive generic constant $\gamma$ depends only on the data if it depends on the parameters $\{N,p_i, C_o, C_1\}$; for the summation notation we refer to Section \ref{Preliminaries}.
\vskip0.1cm \noindent 
Our main results concern the integral inequalities which, for the sake of simplicity, we state in a forward cylinder centered at the origin. 

\vskip0.1cm \noindent 

\noindent First, we state the Harnack-type inequalities for the $L^1_{loc}(\Omega)$ norm of the solution evolving in time, sorting out the case of anisotropic intrinsic geometry from the anisotropic standard one. 

\begin{theorem}[Intrinsic $L^1$-$L^1$ Harnack-type inequality] \label{l1l1}
Let $u$ be a non-negative, local weak solution to equation \eqref{EQ}-\eqref{structure-conditions} in $\Omega_T$, $1<p_i<2$ for all $i=1,\cdots, N$. Let $t, \rho>0$ be such that the inclusion 
\[\K_{2\rho}(t)\times [0,t] \subset \Omega_T,\]
holds true. Then, there exists a positive constant $\gamma$ depending only on the data such that, either there exists an index $i \in \{1, \dots, N\}$ for which
\begin{equation}\label{alternative}
    C^{p_i}\rho^p > \min\{1,\nu^{p-p_i}, \nu^{p}\}, \qquad \text{where} \qquad \nu= \bigg(\frac{t}{\rho^p}\bigg)^{\frac{1}{2-p}},
\end{equation} or, denoting $\lambda= N(p-2)+p$, we have
\begin{equation}\label{l1-l1}
    \sup_{0\leq \tau \leq  t} \int_{\K_{\rho}(t)} u(x,\tau)\, dx \leq \gamma \inf_{0\leq \tau \leq t} \int_{\K_{2\rho}(t)} u(x,\tau) \, dx + \gamma \bigg( \frac{t}{\rho^\lambda} \bigg)^{\frac{1}{2-p}}.
\end{equation} \noindent 
\end{theorem}

\begin{theorem}[Standard $L^1$-$L^1$ Harnack-type inequality] \label{el1l1}
Let $u$ be a non-negative, local weak solution to equation \eqref{EQ}-\eqref{structure-conditions} in $\Omega_T$, $1<p_i<2$ for all $i=1,\cdots, N$. Let $t, \rho>0$ be such that the inclusion 
\[\k_{2\rho}\times [0,t] \subset \Omega_T \] holds true. Then, there exists a positive constant $\gamma$ depending only on the data such that, either there exists an index $i \in \{1, \dots, N\}$ for which
\begin{equation}\label{ealternative}
    C^{p_i}\rho^p > \min\{1,\nu_{\Sigma}^{p_i}\}, \qquad \text{where} \qquad \nu_{\Sigma}= \sum_k \bigg(\frac{t}{\rho^p}\bigg)^{\frac{1}{2-p_k}},
\end{equation} or, denoting $\lambda_i= N(p_i-2)+p$, we have
\begin{equation}\label{el1-l1}
    \sup_{0\leq \tau \leq  t} \int_{\k_{\rho}} u(x,\tau)\, dx \leq \gamma \inf_{0\leq \tau\leq t} \int_{\k_{2\rho}} u(x,\tau) \, dx + \gamma \sum_i \bigg( \frac{t}{\rho^{\lambda_i}} \bigg)^{\frac{1}{2-p_i}}.
\end{equation} 
\end{theorem}

\begin{remark} We remark that in Theorems \ref{l1l1} and \ref{el1l1} the constants $\lambda, \lambda_{i}$ can be of either sign. 
\end{remark}

\noindent Then, considering  extra local regularity assumptions on $u$ such as local boundedness and $u \in L^r_{loc}(\Omega_T)$, for some $r>1$, we have the following $L^r$-$L^{\infty}$ estimates, valid for exponents $p>2N/(N+r)$.

\begin{theorem}[Intrinsic Backwards $L^r$-$L^{\infty}$ estimate] \label{lr-decay} Let $u$ be a non-negative, locally bounded, local weak solution to \eqref{EQ}-\eqref{structure-conditions} in $\Omega_T$, and suppose that for some $r > 1$ it satisfies both $u \in L^r_{loc}(\Omega_T)$ and 
\begin{equation}
\lambda_r= N(p-2)+rp >0.
\end{equation}\noindent Then, there exists a positive constant $\gamma$ depending only on the data, such that for all cylinders 
\[\K_{2\rho}(t) \times [0,t] \subset \Omega_T,\]
either there exists an index $i \in \{1, \dots, N\}$ such that \eqref{alternative} holds true, or 
\begin{equation}\label{decay}
            \sup_{\K_{{\rho}/{2}}(t)\times [t/2,t]} u\, \, 
            \leq \gamma  t^{-\frac{N}{\lambda_r}} \bigg( \int_{\K_{2\rho}(t)} u^r(x,0)\, dx\bigg)^{\frac{p}{\lambda_r}}+ \gamma  \bigg( \frac{t}{\rho^p} \bigg)^{\frac{1}{2-p}} .
    \end{equation}
\end{theorem}

\vskip0.4cm \noindent

\begin{theorem}[Standard Backwards $L^r$-$L^{\infty}$ estimate] \label{elr-decay} Let $u$ be a non-negative, locally bounded, local weak solution to \eqref{EQ}-\eqref{structure-conditions} in $\Omega_T$ and suppose additionally that, for some $r >1$,  $u \in L^r_{loc}(\Omega_T)$ and 
\begin{equation}
\lambda_r= N(p-2)+rp >0.
\end{equation}\noindent Then, there exists a positive constant $\gamma$ depending only on the data, such that for all cylinders 
\[\k_{2\rho} \times [0,t] \subset \Omega_T,\]
either there exists an index $i \in \{1, \dots, N\}$ for which \eqref{ealternative} holds true, or
\begin{equation}\label{edecay}
            \sup_{\k_{{\rho}/{2}}\times [t/2,t]} u\, \, 
            \leq \gamma  t^{-\frac{N}{\lambda_r}} \bigg( \int_{\k_{2\rho}} u^r(x,0)\, dx\bigg)^{\frac{p}{\lambda_r}}+ \gamma  \sum_i \bigg( \frac{t}{\rho^p} \bigg)^{\frac{\lambda_{i,r}}{(2-p_i)\lambda_r}} +  \gamma   \sum_i\bigg(\frac{t}{\rho^p} \bigg)^{\frac{1}{2-p_i}} ,
    \end{equation}\noindent for exponents $\lambda_{i,r}= N(p_i-2)+pr$.
\end{theorem}

\begin{remark} In the prototype degenerate case ($p_i>2$ for all $i=1, \dots, N$) estimates \eqref{decay}-\eqref{edecay} hold true without the second term (and third) on the right-hand side of the inequality (see for instance \cite{CiaGua-Liouville} and \cite{DMV}). Similarly, to what discussed in \cite{DBH-sing}, the distinction between the two approaches relies in the consideration of solutions that are either local or global in time. With the integral Harnack estimates derived in this paper, it is possible to embark on the path of global existence of solutions to \eqref{Prototype}. To this aim we observe that the first term on the right hand side of \eqref{decay} is formally the same as in the degenerate case, while the second term on the right-hand side controls the growth of the solution for large times.
\end{remark}

\noindent Finally, we state the main results or our analysis: Harnack-type estimates considered in the topologies $L^{\infty}_{loc}(\Omega)$ to $L^{1}_{loc}(\Omega)$, again distinguishing when the anisotropic geometry considered is intrinsic or standard.

\begin{theorem}[Intrinsic $L^1$-$L^{\infty}$ Harnack-type inequality] \label{l1-linfty-THM} Let $u$ be a non-negative, locally bounded, local weak solution to \eqref{EQ}-\eqref{structure-conditions} and suppose $p$ is in the supercritical range, i.e. 
\[\lambda =N(p-2)+p>0.\] Then, there exists a positive constant $\gamma$ depending only on the data such that, for all cylinders 
\[\K_{2\rho}(t) \times [0,t] \subset \Omega_T,\]
either there exists $i \in \{1, \dots, N\}$ for which \eqref{alternative} holds true, or 
\begin{equation}\label{l1-linfty-Harnack}
    \sup_{\K_{{\rho}/{2}}(t) \times [t/2\, , \, t]} u \leq \gamma\,  t^{\frac{-N}{\lambda}} \bigg( \inf_{0\leq \tau \leq t}  \int_{\K_{2\rho(t)}} u(x,\tau)\, dx\bigg)^{\frac{p}{\lambda}}+ \gamma \, \bigg( \frac{t}{\rho^p} \bigg)^{\frac{1}{2-p}}.
\end{equation}    
\end{theorem}

\vskip0.4cm \noindent

\begin{theorem}[Standard $L^1$-$L^{\infty}$ Harnack-type inequality] \label{el1-linfty-THM} Let $u$ be a non-negative, locally bounded, local weak solution to \eqref{EQ}-\eqref{structure-conditions} and suppose $p$ is in the supercritical range, i.e. 
\[\lambda =N(p-2)+p>0.\] Then, there exists a positive constant $\gamma$ depending only on the data such that, for all cylinders 
\[\k_{2\rho} \times [0,t] \subset \Omega_T,\]
either there exists $i \in \{1, \dots, N\}$ for which \eqref{ealternative} holds true, or 
\begin{equation}\label{el1-linfty-Harnack}
    \sup_{\k_{{\rho}/{2}} \times [t/2\, , \, t]} u \leq \gamma\,  t^{\frac{-N}{\lambda}} \bigg( \inf_{0\leq \tau \leq t}  \int_{\k_{2\rho}} u(x,\tau)\, dx\bigg)^{\frac{p}{\lambda}}+ \gamma \sum_i \, \bigg( \frac{t}{\rho^p} \bigg)^{\frac{\lambda_i}{(2-p_i)\lambda}}+ \gamma   \sum_i\bigg(\frac{t}{\rho^p} \bigg)^{\frac{1}{2-p_i}},
\end{equation}    \noindent for $\lambda_i=N(p_i-2)+p$.
\end{theorem}

\vskip0.4cm \noindent 

 \subsection*{Rates of Extinction}
\noindent The fact that certain solutions to \eqref{EQ}-\eqref{structure-conditions} with $C=0$ are subject to extinction in finite time has been studied in \cite{AS2} and also in \cite{AS-book2} (we refer to \cite{BenCra},\cite{Diaz-anulacion}, \cite{DB}, for the isotropic case, all  $p_i \equiv p$). In \cite{AS2}, the authors suppose $u$ to be a solution to 
\begin{equation}\label{extinction-system}
\begin{cases}
    \partial_t u - \sum_i \partial_i (a_i(x,t, u) |\partial_i u |^{p_i-2} \partial_i u)=0, & (x,t) \in \Omega \times (0,T),\\
    u=0 &  (x,t)\in \partial\Omega \times (0,T), \\
    u(x,0)=u_0(x) & x\in \Omega,
\end{cases}
\end{equation} with $u_0 \in L^2(\Omega)$ and where $a_i : \Omega \times (0,T)\times \R \rightarrow \R$ are Caratheodory functions satisfying $a_0 \leq a_i(x,t, s) \leq A_0$, for $a_0,A_0>0$ structural constants. Within this framework, the authors show that if $1<p<2$, being $p= N/(\sum_i {p_i}^{-1})$ the harmonic average of the exponents $p_i$, then the energy solutions to \eqref{extinction-system} vanish in a finite time, i.e
\[u(x,t)\equiv 0 \quad \text{for all} \quad t\ge T^*= \bigg( C_e  \|u_0\|_{2, \Omega}^2\bigg)^{\frac{2}{2-p}}, \qquad \qquad C_e=C_e(a_0,A_0,p_i,N)>0.\]

\vskip0.2cm 

\noindent By using a weaker definition of solution (see Definition \ref{def-sol}), here we assume $u$ is a non-negative, local weak solution to \eqref{EQ}-\eqref{structure-conditions} in $\Omega_T$, with $C=0$, $1<p_i<2$ for all $i=1, \dots, N$, and that there exists an extinction time $T^*<T$ for $u$. Then, similarly to \cite{DBKW}, we use the $L^1$-$L^1$ Harnack-type inequalities \eqref{l1-l1}-\eqref{el1-l1} to evaluate the decay of the $L^1_{loc}(\Omega)$ norm of $u$ toward its extinction and the $L^1$-$L^{\infty}$ Harnack-type inequalities \eqref{l1-linfty-Harnack}-\eqref{el1-linfty-Harnack} to estimate the rate of extinction of the solution in a whole half cylinder approaching $T^*$. These two properties require different assumptions on the exponents $p_i$. We divide the cases distinguishing the underlying geometry.

\subsubsection*{Intrinsic Geometry} \vskip0.1cm \noindent  Let $\tau,\rho>0$ be fixed such that $\K_{4\rho}(T^*-\tau) \subseteq \Omega$. 

\begin{itemize}
    \item The mass decays within the law
    \[\|u(\cdot, \tau)\|_{1,(\K_{\rho}(T^*-\tau))}= \int_{\K_{2\rho}(T^*-\tau)} u(x,\tau) \, dx \leq \gamma \bigg(\frac{T^*-\tau}{\rho^{\lambda}}  \bigg)^{\frac{1}{2-p}}, \]
    for a positive constant $\gamma$ depending only on the data. Hence the mass $\|u(\cdot, \tau)\|_{L^1(\K_{\rho}(T^*-\tau))}$ of the solution locally decays (to zero) as $(T^*-\tau)^{1/(2-p)}$ in a space configuration depending on time but with unchanged measure $|\K_{\rho}(T^*-\tau)|=(2\rho)^N$. 
\vskip0.2cm \noindent 

    \item If $\lambda=N(p-2)+p>0$, then the solution has the following vanishing rate:
    \[ \sup_{\K_{\rho}(T^*-\tau)\times [(T^*+\tau)/2,\,T^{*}]} u \leq \gamma \bigg( \frac{T^*-\tau}{\rho^p} \bigg)^{\frac{1}{2-p}}, \quad \forall \tau \in (0,\,T^*),\]
    for a positive constant $\gamma$ depending only on the data. Choosing $T^*/2<t<T^*$, it is possible to specialize this decay to an ultra-contractive bound
    \[\|u(\cdot, t)\|_{\infty, \K_{\rho}(T^*-t)} \leq \gamma \bigg( \frac{T^*-t}{\rho^p} \bigg)^{\frac{1}{2-p}}. \]
    \noindent This estimate shows that the rate of local decay of the $L^{\infty}$-norm of the solution, in a space configuration depending on each time $t$, is again of the type $(T^*-t)^{1/(2-p)}$ but now for a different power of the radius $\rho$.  
\end{itemize}

\vskip0.2cm 
\noindent We observe that when $t\rightarrow T^*$ the time intrinsic cube $\K_{\rho}(T^*-t)$ shrinks along the directions $x_k$ for which $p_k>p$, while in the other directions it stretches to infinity; this particular phenomenon occurs keeping the measure $|\K_{\rho}(T^*-t)|$ unchanged. Therefore, the inclusion $\K_{4\rho}(T^*-t)\subseteq \Omega$ degenerates according to the choice of time.

\vskip0.2cm

\subsubsection*{Standard Anisotropic Geometry}
For a positive number $\rho$, let us consider the anisotropic standard cube $\k_{\rho}$ as in \eqref{anisociccio}, for $\rho>0$ such that $\k_{\rho} \subset \Omega$. We can estimate the local decay of its $L^1$ and $L^\infty$ norms as above, but this time in a space geometry that is time independent, paying the price of having more involved estimates.
\vskip0.2cm \noindent 
\begin{itemize}
    \item Description of the mass decay 
\[ \|u(\cdot, \tau)\|_{L^1(\k_{\rho})}=\int_{\k_\rho} u(x,\tau) \, dx  \leq \gamma \sum_{i}\bigg(\frac{T^*-\tau}{\rho^{\lambda_i}}\bigg)^{\frac{1}{2-p_i}}, \quad \forall \quad 0<\tau \leq T^*.\]
\noindent When considering times $\tau$ approaching $T^*$, the mass of the solution $\|u(\cdot, \tau)\|_{L^1(\k_{\rho})}$  decays to zero at the rate $(T^*-\tau)^{1/(2-p_N)}$, while when considering larger times $(T^*-\tau)>1$ the rate is $(T^*-\tau)^{1/(2-p_1)}$. 

\item For any time $0< \tau<T^*$, and assuming that  $\lambda>0$, we have a description of the local decay of the essential supremum of the solution as
\[\sup_{\k_{\rho} \times [(T^* +\tau)/2 , \, T^*]} u  \leq \gamma \sum_i \, \bigg( \frac{T^* -\tau}{\rho^p} \bigg)^{\frac{\lambda_i}{(2-p_i)\lambda}} + \gamma \sum_i \, \bigg( \frac{T^* -\tau}{\rho^p} \bigg)^{\frac{1}{2-p_i}}, 
  \] for $\gamma$ positive constant depending only on the data $\{C_o,C_1,C_2,  p_i,N\}$  and being $\lambda_i=N(p_i-2)+p$. Here we observe that a decay rate towards extinction, i.e. for times $(T^*-\tau)<1$, is given from this estimate only with the extra assumption $\lambda_i=N(p_i-2)+p>0$ for all $i=1, \dots, N$, and the solution vanishes in the half-cylinder as fast as $(T^*-\tau)^{\lambda_1/[(2-p_N)\lambda]}$. This behavior is confirmed by those solutions that are constant along $N-1$ space coordinates and behave like a $p_1$ or $p_N$-Laplacian by means of the only free variable.\end{itemize}

\section*{Acknowledgements}
\noindent The third author is partially supported by the Grant EFDS-FL2-08 of the found The European Federation of Academies of Sciences and Humanities (ALLEA) and
by the Project ”Mathematical modelling of complex dynamical systems and processes caused
by the state security” (Reg. No. 0123U100853). The second author  was financed by Portuguese Funds through FCT - Funda\c c\~ao para a Ci\^encia e a Tecnologia - within the Projects UIDB/00013/2020 and UIDP/00013/2020. The first author acknowledges the support of the department of Mathematics of the University of Bologna Alma Mater and the Italian PNR (MIUR) fundings 2021-2027.

\section{Functional Setting and Notation} \label{Preliminaries}\noindent
\subsection*{Functional Setting} \noindent We define the anisotropic spaces of locally integrable functions as
\[W^{1,{\bf p}}_{loc}(\Omega)= \{u \in W^{1,1}_{loc}(\Omega)\, |\, \, \partial_i u \in L^{p_i}_{loc}(\Omega) \},\]
\[L^{{\bf p}}_{loc}(0,T; W^{1, {\bf p}}_{loc}(\Omega))= \{u \in L^1_{loc}(0,T; W^{1,1}_{loc}(\Omega)\, |\, \, \partial_i u \in L_{loc}^{p_i}(0,T; L^{p_i}_{loc} (\Omega))\},\]
and the respective spaces of functions with zero boundary data
\[W^{1,{\bf p}}_{o}(\Omega)= \{u \in W^{1,1}_{o}(\Omega)\, |\, \, \partial_i u \in L^{p_i}_{loc}(\Omega) \},\]
\[L^{{\bf p}}_{loc}(0,T; W^{1, {\bf p}}_{o}(\Omega))= \{u \in L^1_{loc}(0,T; W^{1,1}_{o}(\Omega)\, |\, \, \partial_i u \in L_{loc}^{p_i}(0,T; L^{p_i}_{loc} (\Omega))\}.\]
\noindent It is known (see  \cite{IlinBesovNikolski}, \cite{russian1}) that when $p>N$ the embedding  $W^{1,{\bf p}}(\Omega) \hookrightarrow C^{0,\alpha}_{loc}(\Omega)$ for $\Omega$ regular enough. Therefore in this work we will consider $p<N$.
\begin{definition} \label{def-sol}
A function 
\[u \in C(0,T; L^2_{loc}(\Omega))\cap L^{\bf p}_{loc}(0,T; W^{1,{\bf p}}_{loc}(\Omega))\] 
 is called a local weak sub(super)-solution to \eqref{EQ} in $\Omega_T$ if, for all times $0\leq t_1 \leq t_2 \leq T$ and for all compact sets $K\subset \subset \Omega$, it satisfies the inequality 
\begin{equation} \label{EQ-weak}
\begin{aligned}
\int_K u \varphi \, dx\, \bigg|_{t_1}^{t_2} + \int_{t_1}^{t_2} \int_K \{-u \partial_{\tau}\varphi + & \sum_i A_i(x,t,u,Du)  \partial_i  \varphi \} dxd\tau\\& 
\leq (\ge) \int_{t_1}^{t_2} \int_K B(x,t,u,Du) \varphi dx\,d\tau,
\end{aligned}
\end{equation} for all non-negative 
test functions $\varphi \in W^{1,2}_{loc}(0,T; L^2_{loc}(\Omega)) \cap L^{\bf p}_{loc} (0,T; W^{1,{\bf p}}_{o} (\Omega))$.
\end{definition}\noindent 
This last membership of the test functions, together with the structure conditions \eqref{structure-conditions}, ensure that all the integrals in \eqref{EQ-weak} are finite. Moreover, as $\varphi$ vanishes along the lateral boundary of $\Omega_T$, its integrability increases thanks to the following known embedding theorem.

\begin{lemma}(Anisotropic Gagliardo-Sobolev-Nirenberg, \cite{DMV})
\label{embedding} \vskip0.1cm \noindent 
Let $\Omega\subseteq \R^N$ be a rectangular domain, $p<N$, and $\sigma\in [1, p^{*}]$. For any number $\theta\in [0,\,  p/  p^{*}]$ define 
\[
q=q(\theta, {\bf p})=\theta \,  p^{*}+\sigma\, (1-\theta),
\]
Then there exists a positive constant $c=c(N, {\bf p}, \theta, \sigma)>0$ such that
\begin{equation}
\label{PS}
\iint_{\Omega_{T}}|\varphi|^{q}\, dx\, dt\leq c\, T^{1-\theta\, \frac{ p^{*}}{ p}}\left(\sup_{t\in (0, T]}\int_{\Omega}|\varphi|^{\sigma}(x, t)\, dx\right)^{1-\theta}\prod_i \left(\iint_{\Omega_{T}}|\partial_i \varphi|^{p_{i}}\, dx\, dt\right)^{\frac{\theta\,   p^{*}}{N \, p_{i}}},
\end{equation} 
for any $\varphi\in L^{1}(0, T; W^{1,1}_{o}(\Omega))$, being the inequality trivial when the right-hand side is unbounded.
\end{lemma} \noindent 

\subsection*{Notation} In what follows we introduce the notation we will be using along the text.

\vskip0.2cm 
\noindent 

\begin{itemize}
\item[] We shorten the notation on sums and products when they are intended for all indexes $i,j,k \in \{1, \dots, N\}$, 
\[\sum_i:= \sum_{i=1}^N\quad  \qquad \text{and} \qquad \quad \prod_i:= \prod_{i=1}^N  \quad. \]

\noindent Only when the sum runs over a different range of exponents will be further specified.

\vskip0.2cm \noindent 
\item[] Exponents are ordered, \[1<p_1 \leq p_2 \leq \dots, \leq p_N<2,\]
and $p$ stands for the harmonic average

  \[\displaystyle{p:=\bar{p}=  N/\bigg(\sum_i 1/p_i}\bigg)  \, \, .\]

    \item[] We denote by $\partial_i u$ the weak directional space derivatives and by $\partial_t u$ the weak time-derivative (see \eqref{Steklov} for more details). Finally, $\nabla u= (\partial_1u, \dots, \partial_N u)$.
\vskip0.3cm \noindent 

    \item[] Our geometrical setting will distinguish between two types of $N$-dimensional cubes:

    \begin{itemize}
        \item[] Anisotropic intrinsic cube \[
    \K_{a \rho}(t):=\prod_i \bigg{\{} |x_i|< a \ \rho^{\frac{p}{p_i} \frac{(2-p_i)}{(2-p)}} t^{\frac{(p_i-p)}{(2-p)p_i}} \bigg{\}}, \quad a>0 \ , \qquad |\K_{\rho}(t)|=(2\rho)^N 
    \]
        \item[] Anisotropic standard cube 
    \[
    \k_{a \rho}:=\prod_i \bigg{\{} |x_i|< (a  \rho)^{\frac{p}{p_i}} \bigg{\}}, \quad a>0 \ , \qquad |\k_{\rho}|=(2\rho)^N.
    \]
    \end{itemize}

    \item[] We will use two exponents for the decay rates:
    \[\lambda_r= N(p-2)+rp     \qquad \text{\&} \qquad \lambda_{i,r}= N(p_i-2)+rp,\]
    when $r=1$, the subscript $r$ is dropped writing $\lambda= N(p-2)+p$ and $\lambda_{i}=N(p_i-2)+p.$
    \vskip0.3cm \noindent 
    
\item[] Given a measurable function $u :E \subset \R^{N+1} \rightarrow \R$, we denote by $\sup_E u$ ($\inf_E u$) the essential supremum (essential infimum of $u$) in $E$ with respect to the Lebesgue measure.

\vskip0.3cm \noindent
\item[] We denote by $\gamma$ a generic positive constant that depends only on the structural data $\{p_i, N, C_o, C_1\}$ to  \eqref{EQ}-\eqref{structure-conditions}, and it may vary in the estimate from line to line.

\vskip0.3cm \noindent
\item[] Young's Inequality Convention. In our estimates we will repeatedly use Young's inequality in the following form: for $q>1$ and $a,b, \epsilon>0$ fixed, we use the well-known inequality
\begin{equation} \label{Young}
    ab \leq \epsilon a^q+ \gamma(\epsilon) b^{q'} \ , 
\end{equation}
\[ \text{with}  \qquad  \qquad q'=(1-1/q)^{-1}, \quad \text{and} \quad \gamma(\epsilon)= \bigg( \frac{q-1}{q^{1/(q-1)}q}\bigg)  \bigg( \frac{1}{\epsilon} \bigg)^{\frac{1}{q-1}}.\]
\noindent The constant $\epsilon$ will not be specified as long as it depends only on the data $\{p_i, N, C_o, C_1\}$.
\end{itemize}

\section{Proof of $L^1$-$L^1$ Harnack estimates} \label{l1l1proofs}

\noindent In this Section we prove Theorems \ref{l1l1}-\ref{el1l1}, dividing the argument whether the anisotropic space geometry considered is the standard or the intrinsic one.

   \subsection*{Intrinsic Anisotropic Geometry: Proof of Theorem \ref{l1l1}}

\noindent We consider a fixed time-length $0<t<T$, and let $\rho>0$ be small enough to allow the inclusion  
\begin{equation*}\label{anisocilindro}
\begin{aligned}
\Q_{\rho}(t):=\K_{\rho}(t)\times [0, t] &=\prod_i \bigg{\{} |x_i|< \rho^{\frac{p}{p_i}} \nu^{\frac{p_i-p}{p_i}} \bigg{\}}\times [0, t]  \subseteq \Omega_T,
\end{aligned} \end{equation*}
\noindent for the fixed quantity
\begin{equation}\label{nu}
\nu= \bigg(\frac{t}{\rho^p}  \bigg)^{\frac{1}{2-p}}.
\end{equation} \noindent

\begin{lemma}
Let $u$ be a non-negative local weak super-solution to \eqref{EQ} in $\Omega_T$ and $\sigma \in (0,1)$ a number. Then, there exists a positive constant $\gamma$ depending only on the data such that, either \eqref{alternative} holds true for some $i =1, \dots, N$, or we have 
\begin{equation}\label{estimate1}
\sum_i \frac{1}{\rho^{\frac{p}{p_i}}} \bigg( \frac{t}{\rho^p}\bigg)^{\frac{p-p_i}{p_i(2-p)}}\int_{0}^{t} \int_{\K_{\sigma \rho}(t)} |\partial_i u |^{p_i-1}\, dx d\tau \leq \frac{\gamma}{(1-\sigma)^p} \sum_i \bigg(\frac{t}{\rho^{\lambda}} \bigg)^{\frac{2-p_i}{p_i(2-p)}} \bigg{\{} S+ \bigg(\frac{t}{\rho^{\lambda}}  \bigg)^{\frac{1}{2-p}}  \bigg{\}}^{\frac{2(p_i-1)}{p_i}},
\end{equation} being \begin{equation*}
S=\sup_{0\leq \tau \leq t}  \int_{\K_{\rho}(t)} u(x,\tau) \ dx \quad \text{and} \quad \lambda= N(p-2)+p.\\
\end{equation*}
\end{lemma}

\begin{proof} For each $i= 1, \dots, N$ we apply H\"older's inequality to the quantity to be estimated, 
\begin{equation*}\label{est1}
\begin{aligned}
& \int_{0}^{t}  \int_{\K_{\sigma \rho}(t)} |\partial_i u |^{p_i-1}\, dx d\tau \\
&=  \int_{0}^{t} \int_{\K_{\sigma \rho}(t)} \bigg( |\partial_i u |^{p_i-1} \tau^{\frac{1}{p_i} (\frac{p_i-1}{p_i})} (u+ \nu)^{\frac{-2}{p_i}(\frac{p_i-1}{p_i})} \bigg) \, \bigg( \tau^{\frac{-1}{p_i} (\frac{p_i-1}{p_i})} (u+ \nu)^{\frac{+2}{p_i}(\frac{p_i-1}{p_i})} \bigg) \, dx d\tau \\
&\leq   \bigg(  \int_{0}^{t} \int_{\K_{\sigma \rho}(t)} |\partial_i u |^{p_i} \tau^{\frac{1}{p_i}} (u+\nu)^{\frac{-2}{p_i}}\, dxd\tau \bigg)^{\frac{p_i-1}{p_i}} \bigg(   \int_{0}^{t} \int_{\K_{\sigma \rho}(t)}  \tau^{\frac{-1}{p_i}(p_i-1)} (u+\nu)^{\frac{2}{p_i}(p_i-1)} \, dx d\tau \bigg)^{\frac{1}{p_i}}  \\
& =:  I_{1,i}^{\frac{p_i-1}{p_i}} I_{2,i}^{\frac{1}{p_i}}\,. 
\end{aligned}
\end{equation*}
\noindent Next, we estimate $I_{2,i}$ by taking the supremum in time and then using H\"older's inequality
\begin{equation*}\label{est-I2}
\begin{aligned}
    I_{2,i} & = \int_{0}^{t} \int_{\K_{\sigma \rho}(t)}  \tau^{\frac{1}{p_i}-1} (u+\nu)^{\frac{2}{p_i}(p_i-1)} \, dx d\tau \\
    &\leq \int_{0}^{t} \tau^{\frac{1}{p_i}-1} d\tau \bigg( \sup_{0 \leq \tau \leq t} \int_{\K_{\rho}(t)} (u(\tau)+\nu)^{\frac{2}{p_i}  (p_i-1)}\, dx \bigg) \\
    & \leq \gamma\,  t^{\frac{1}{p_i}} \, |\K_{\rho}(t)|^{\frac{2-p_i}{p_i}}\, \bigg( \sup_{0\leq \tau \leq  t}  \int_{\K_{\rho}(t)} (u(\tau)+\nu) \, dx \bigg)^{\frac{2(p_i-1)}{p_i}}   
    \\
    &= \gamma \, t^{\frac{1}{p_i}} \, \rho^{N\frac{(2-p_i)}{p_i}}\, \bigg( \sup_{0\leq \tau \leq  t}  \int_{\K_{\rho}(t)} u(\tau) \, dx+ \nu \rho^N \bigg)^{\frac{2(p_i-1)}{p_i}} 
    \\
    &=: \gamma \, \bigg( \frac{t}{\rho^{\lambda_i}} \bigg)^{\frac{1}{p_i}} \, \rho^{\frac{p}{p_i}} \bigg{\{} S+\bigg( \frac{t}{\rho^{\lambda}}  \bigg)^{\frac{1}{2-p}}\bigg{\}}^{\frac{2(p_i-1)}{p_i}}.
\end{aligned}
\end{equation*}
In the last steps we have used the property $|\K_{\rho}(t)|=(2\rho)^N$ and the definition of $\nu$, $\lambda_i,\lambda$ (see the statement of Theorem \ref{l1l1}). Now we estimate $I_{1,i}$ using the inequalities \eqref{Negative-EE} within the considered geometry: we test indeed repeatedly, for $i=1,\dots,N$, 
 equation \eqref{EQ} with the function 
\begin{equation*} \varphi_i(x,\tau)= -\tau^{\frac{1}{p_i}} (u(x,\tau)+\nu)^{1-\frac{2}{p_i}} \zeta(x), \qquad \zeta(x)= \prod_i \zeta_i(x_i)^{p_i},\quad \hat{\zeta}^{j}:= \prod_{i\ne j} \zeta_i(x_i)^{p_i}\end{equation*}\noindent being $\zeta$ a smooth cut-off function between the sets $\K_{\sigma \rho}(t)$ and $\K_{\rho}(t)$, hence enjoying the properties
\begin{equation} \label{test-radii} 0\leq \zeta\leq 1, \qquad ||\partial_i \zeta_i||_{\infty}\leq \gamma \bigg( [(1-\sigma)\rho]^{\frac{p}{p_i}} (t/\rho^p)^\frac{(p_i-p)}{(2-p)p_i} \bigg)^{-1}= \gamma/\bigg([(1-\sigma)\rho]^{\frac{p}{p_i}} \nu^{\frac{(p_i-p)}{p_i}}\bigg).
\end{equation} \noindent  The number $\nu \in \R^+$ is fixed, and by implementing \eqref{test-radii} into \eqref{Negative-EE} we obtain
\begin{equation}\label{o}
    \begin{aligned}
        \int_{0}^{t} \int_{\K_{\rho}(t)} \sum_j& |\partial_j u |^{p_j}\tau^{\frac{1}{p_i}} (u+\nu)^{-\frac{2}{p_i}} \zeta\, dxd\tau \leq \gamma t^{\frac{1}{p_i}} \int_{\K_{\rho}(t)} (u+\nu)^{\frac{2(p_i-1)}{p_i}}\zeta\, dx\\
        & +   \gamma \sum_j  \frac{\nu^{p-p_j}}{[(1-\sigma)\rho]^p}\bigg[1+\bigg( \frac{C^{p_j} \rho^{{p}}}{\nu^{{p-p_j}}}\bigg)
        \bigg] \int_{0}^{t} \int_{\K_{\rho}(t)}  (u+\nu)^{p_j-\frac{2}{p_i}} \tau^{\frac{1}{p_i}}\, dxd\tau \\
        & + \gamma \sum_j C^{p_j} \int_{0}^{t} \int_{\K_{\rho}(t)} (u+\nu)^{-\frac{2}{p_i}} \tau^{\frac{1}{p_i}}\, dxd\tau =: I_1+I_2+I_3.
    \end{aligned}
\end{equation}

\noindent \noindent Now we manipulate the terms of \eqref{o}, with the aim of obtaining an homogeneous estimate similar to $I_{2,i}$.\vskip0.2cm 
 \noindent The first term on the right is bounded from above by a similar estimate as the one for $I_{2,i}$. \vskip0.2cm \noindent 

\noindent The second term is the one most related with our anisotropic problem; it is here that we specialize our estimates toward homogeneity. We dominate it from above by using $p_i<2$, with the usual trick \[(u+\nu)^{p_j-\frac{2}{p_i}}=(u+\nu)^{\frac{2(p_i-1)}{p_i}}(u+\nu)^{p_j-2} \leq (u+\nu)^{\frac{2(p_i-1)}{p_i}}\nu^{p_j-2},\]  in order to give an homogeneous estimate with respect to $j$-th index, namely
\begin{equation*} \label{crucial}
    \begin{aligned}
       I_2  & = \gamma\sum_j  \frac{\nu^{p-p_j}}{[(1-\sigma)\rho]^p}\bigg[1+\bigg( \frac{C^{p_j} \rho^{{p}}}{\nu^{{p-p_j}}}\bigg)
        \bigg]  \int_{0}^{t} \int_{\K_{\rho}(t)}  (u+\nu)^{p_j-\frac{2}{p_i}} \tau^{\frac{1}{p_i}}\, dxd\tau\\
        & \leq \gamma \sum_j  \bigg[1+ \bigg( \frac{C^{p_j} \rho^{{p}}}{\nu^{{p-p_j}}}\bigg)
        \bigg]\frac{\nu^{p_j-2}}{[(1-\sigma)\rho]^p \nu^{p_j-p}} t^{1+\frac{1}{p_i}} \bigg( \sup_{0\leq \tau \leq t} \int_{K_{\rho}(t)}  (u+\nu)^{\frac{2(p_i-1)}{p_i}}\, dx\bigg)  \\
        & \leq N \gamma \bigg[1+\sum_j \bigg( \frac{C^{p_j} \rho^{{p}}}{\nu^{{p-p_j}}}\bigg)
        \bigg]\bigg(  \frac{t\nu^{p-2}}{\rho^p}\bigg) \bigg(\frac{t}{\rho^{\lambda_i}} \bigg)^{\frac{1}{p_i}} \rho^{\frac{p}{p_i}} \bigg{\{}S+ \nu \rho^N  \bigg{\}}^{\frac{2(p_i-1)}{p_i}},
        \end{aligned}
\end{equation*}
where $\lambda_i= N(p_i-2)+p$, for $i=1,\dots,N$.

\noindent Referring again to \eqref{o}, each $j$-th term of $I_3$ on the right can be estimated by  
\begin{equation*}\label{bacalao} \begin{aligned} 
C^{p_j}\int_0^t  \int_{\K_{\rho}(t)} &(u+\nu)^{\frac{-2}{p_i}} \tau^{\frac{1}{p_i}}\, dxd\tau \\
& \leq  C^{p_j}t^{1+\frac{1}{p_i}}\nu^{-2} \bigg( \sup_{0\leq \tau \leq t}\int_{\K_{\rho}(t)} (u + \nu)^{\frac{2(p_i-1)}{p_i}} \, dx\bigg)\\
&\leq  \bigg( \frac{C^{p_j}\rho^p}{\nu^p} \bigg) \bigg(  \frac{t \nu^{p-2}}{\rho^p}\bigg) \bigg(\frac{t}{\rho^{\lambda_i}} \bigg)^{\frac{1}{p_i}} \rho^{\frac{p}{p_i}} \bigg{\{}S+ \nu \rho^N  \bigg{\}}^{\frac{2(p_i-1)}{p_i}},
\end{aligned} \end{equation*} where the first inequality uses $(u+\nu)^{-2} \leq \nu^{-2}$ and the last inequality is brought similarly to the one for $I_{2,i}$. Finally, collecting everything together we arrive, for each $i=1, \dots, N$, to the estimate 
\begin{equation*}
    \begin{aligned}
    & I_{1,i}
    \leq \frac{\gamma \rho^{\frac{p}{p_i}}}{(1-\sigma)^p} \bigg{\{} 1+\bigg[1+\sum_j \bigg( \frac{C^{p_j}\rho^p}{\nu^{p}} \bigg)+ \bigg(\sum_{j}\frac{C^{p_j}\rho^p}{\nu^{p-p_j}}  \bigg)  \bigg] \bigg( \frac{t\nu^{p-2}}{\rho^p} \bigg) \bigg{\}} \bigg( \frac{t}{\rho^{\lambda_i}} \bigg)^{\frac{1}{p_i}} \bigg{\{} S+\nu \rho^N \bigg{\}}^{\frac{2(p_i-1)}{p_i}}.
    \end{aligned}
\end{equation*}
\noindent If condition \eqref{alternative} is violated for all $i=1, \dots, N$, then the term in squared brackets on the right-hand side is smaller than 3, recalling \eqref{nu}. Thence we go back to the initial estimate and evaluate
\begin{equation*} \label{est-F}
\begin{aligned}
& \int_{0}^{t}  \int_{\K_{\sigma \rho}(t)} |\partial_i u |^{p_i-1}\, dx d\tau \leq  I_{1,i}^{\frac{p_i-1}{p_i}} I_{2,i}^{\frac{1}{p_i}}\\
&\leq \gamma    \bigg(\frac{\rho^{\frac{p}{p_i}}}{(1-\sigma)^p}\bigg( \frac{t}{\rho^{\lambda_i}} \bigg)^{\frac{1}{p_i}} \bigg{\{} S+\nu \rho^N \bigg{\}}^{\frac{2(p_i-1)}{p_i}}  \bigg)^{\frac{p_i-1}{p_i}} \bigg( \, \bigg( \frac{t}{\rho^{\lambda_i}} \bigg)^{\frac{1}{p_i}} \, \rho^{\frac{p}{p_i}} \bigg{\{} S+\bigg( \frac{t}{\rho^{\lambda}}  \bigg)^{\frac{1}{2-p}}\bigg{\}}^{\frac{2(p_i-1)}{p_i}}\bigg)^{\frac{1}{p_i}}\\
&\leq  \gamma  \frac{\rho^{\frac{p}{p_i}} }{(1-\sigma)^p}\bigg( \frac{t}{\rho^{\lambda_i}} \bigg)^{\frac{1}{p_i}} \bigg{\{} 
S+\bigg( \frac{t}{\rho^\lambda} \bigg)^{\frac{1}{2-p}}\bigg{\}}^{\frac{2(p_i-1)}{p_i}}
\end{aligned}
\end{equation*}
and thereby
\[ \sum_i \frac{1}{\rho^{\frac{p}{p_i}}} \bigg( \frac{t}{\rho^p}\bigg)^{\frac{p-p_i}{p_i(2-p)}}\int_{0}^{t}  \int_{\K_{\sigma \rho}(t)} |\partial_i u |^{p_i-1}\, dx d\tau 
\leq   \frac{\gamma}{(1-\sigma)^p} \sum_i \bigg( \frac{t}{\rho^{\lambda}} \bigg)^{\frac{2-p_i}{p_i (2-p)}}   \bigg{\{} 
S+\bigg( \frac{t}{\rho^\lambda} \bigg)^{\frac{1}{2-p}}\bigg{\}}^{\frac{2(p_i-1)}{p_i}} \ .\]

\end{proof}


  \subsection*{Proof of Theorem \ref{l1l1} Concluded}
  Now we perform an iteration on $\sigma \in (0,1)$: we define the increasing radii
  \[\rho_{n,i} := \rho^{\frac{p}{p_i}} \bigg( \frac{t}{\rho^p} \bigg)^{\frac{(p_i-p)}{(2-p)p_i}} \bigg( \sum_{k=0}^n 2^{-k}\bigg),\qquad \qquad \rho_{n+1,i}-\rho_{n,i}=2^{-(n+1)}\,  \rho^{\frac{p}{p_i}} \bigg( \frac{t}{\rho^p} \bigg)^{\frac{(p_i-p)}{(2-p)p_i}},\]
  and consider the family of concentric intrinsic anisotropic cubes
  \[
  \K_n= \prod_i \bigg{\{}|x_i|<  \rho_{n,i}\bigg{\}}, \qquad \tilde{\K}_n=  \prod_i \bigg{\{}|x_i|<  \frac{\rho_{n+1,i}+\rho_{n,i}}{2}\bigg{\}},\qquad \text{with} \]
  \[\qquad  \quad  \K_{\rho}(t) = \K_0 \subset \K_n \subset \tilde{\K}_n \subset \K_{n+1} \subset 
 \K_{\infty} = \K_{2 \rho}(t)= \prod_i \bigg{\{} |x_i|< 2 \rho^{\frac{p}{p_i}} (t/\rho^p)^{\frac{(p_i-p)}{(2-p)p_i}}\bigg{\}}.
  \]
For every $n \in \N\cup \{0\}$, consider time-independent cut-off functions $\zeta_n$ as in \eqref{zeta} between $\K_n$ and $\tilde{\K}_n$, hence satisfying
\[
||\partial_i \zeta_n||_{\infty} \leq \frac{\gamma}{|\rho_{n+1,i}-\rho_{n,i}|} \leq \gamma 2^{n+1}\nu^{\frac{p-p_i}{p_i}}/\rho^{\frac{p}{p_i}}.
\]
We test equation \eqref{EQ} with $\zeta_n$ and we integrate over $\tilde{\K}_n \times [\tau_1, \tau_2]$, for arbitrary time levels $0\leq \tau_1 < \tau_2 \leq t$, to get
\begin{equation}\label{start}
\begin{aligned}
    \int_{\tilde{\K}_n} u(x, \tau_1)\, dx \leq& \int_{\tilde{\K}_{n}} u(x,\tau_2)\, dx \\
    & +\gamma 2^{n+1}\sum_i \bigg(\frac{\nu^{{(p-p_i)}}}{\rho^{p}} \bigg)^{\frac{1}{p_i}} \bigg(C_1+ \bigg(\frac{C^{p_i}\rho^{p}}{\nu^{p-p_i}}\bigg)^{\frac{1}{p_i}}\bigg) \int_{\tau_1}^{\tau_2} \int_{\tilde{\K}_n} |\partial_i u|^{p_i-1}\,dx\, d\tau\\
    &+ \gamma \sum_i 2^{n+1} \bigg(  C^{p_i-1}\bigg(\frac{\nu^{p-p_i}}{\rho^{p}}\bigg)^{\frac{1}{p_i}}+C^{p_i}\bigg) \int_{\tau_1}^{\tau_2} \int_{\tilde{\K}_n} \, dxd\tau.
    \end{aligned}
\end{equation} Assume condition \eqref{alternative} is contradicted for all $i\in \{1, \dots, N\}$; then the second term in parenthesis on the right of \eqref{start} is bounded above by $C_1 +1$, while the third term is estimated by
\begin{equation*}
    \begin{aligned}
        \gamma 2^{n+1} \sum_i \bigg( C^{p_i-1}\bigg(\frac{\nu^{p-p_i}}{\rho^{p}}\bigg)^{\frac{1}{p_i}}&+C^{p_i}\bigg) t \rho^N \\
        &= \gamma 2^{n+1}\sum_i \bigg[\bigg( \frac{C^{p_i} \rho^p}{\nu^p}\bigg)^{\frac{(p_i-1)}{p_i}} + \bigg( \frac{C^{\frac{p_i}{p}}\rho}{\nu} \bigg)^{p-1} (C^{\frac{p_i}{p}} \rho)\bigg] \bigg( \frac{t}{\rho^{\lambda}}\bigg)^{\frac{1}{2-p}}\\&
        \leq \gamma 2^{n+1} \bigg( \frac{t}{\rho^{\lambda}} \bigg)^{\frac{1}{2-p}}.
    \end{aligned}
\end{equation*}
\noindent Putting all the pieces together we obtain the estimate 

\begin{equation}\label{restart}
\begin{aligned}
    \int_{\tilde{\K}_n} u(x, \tau_1)\, dx \leq& \int_{\tilde{\K}_{n}} u(x,\tau_2)\, dx\\
    & +\gamma 2^{n}\sum_i \bigg(\frac{\nu^{{(p-p_i)}}}{\rho^{p}} \bigg)^{\frac{1}{p_i}} \int_{\tau_1}^{\tau_2} \int_{\tilde{\K}_n} |\partial_i u|^{p_i-1}\,dx\, d\tau+ \gamma 2^{n} \bigg( \frac{t}{\rho^{\lambda}} \bigg)^{\frac{1}{2-p}}.
    \end{aligned}
\end{equation} 
\noindent By continuity of $u$ as a map $[0,T]\rightarrow L^{2}_{loc}(\Omega)$, we take $\tau_2$ as the time level in $[0,t]$ such that 
\[ \I= \inf_{0\leq \tau \leq t} \int_{2\K_{\rho}(t)} u(x, \tau)\, dx = \int_{2\K_{\rho}(t)} u(x, \tau_2)\, dx,\]
and $\tau_1$ as the time level satisfying
\[
S_n:= \sup_{0\leq \tau \leq t} \int_{\K_n} u(x,\tau)\, dxd\tau=  \int_{\K_n} u(x,\tau_1)\, dxd\tau.
\]
It is precisely for this choice of ordering between $\tau_1$ and $\tau_2$ that we need $u$ to be a {\it solution}, and not only a super-solution. Now we evaluate the second term in \eqref{restart} with the inequality \eqref{estimate1} applied to the pair of cylinders $\tilde{\K}_n \times [0,t] \subset \K_{n+1}\times [0,t]$ and develop the definition of $\nu$ to write 
\begin{equation*} 
    \begin{aligned}
        S_{n} &\leq \I+ \gamma b^{n} \sum_i \nu^{\frac{p-p_i}{p_i}}\bigg(\frac{t}{\rho^{\lambda_i}} \bigg)^{\frac{1}{p_i}} \bigg{\{} S_{n+1}+ \bigg(\frac{t}{\rho^{\lambda}}  \bigg)^{\frac{1}{2-p}}  \bigg{\}}^{\frac{2(p_i-1)}{p_i}}+\gamma  2^{n} \bigg( \frac{t}{\rho^{\lambda}} \bigg)^{\frac{1}{2-p}}\\
        & \leq  \I+ \gamma b^n \sum_i \bigg(\frac{t}{\rho^{\lambda}}  \bigg)^{\frac{(2-p_i)}{(2-p)p_i}} \bigg{\{} S_{n+1}+ \bigg(\frac{t}{\rho^{\lambda}}  \bigg)^{\frac{1}{2-p}}  \bigg{\}}^{\frac{2(p_i-1)}{p_i}}+\gamma  2^{n} \bigg( \frac{t}{\rho^{\lambda}} \bigg)^{\frac{1}{2-p}},\qquad b= 2^{p+1}>1.
    \end{aligned}
\end{equation*} By using Young's inequality on each $i$-th term with exponents $\frac{2(p_i-1)}{p_i}+\frac{2-p_i}{p_i}=1$ we get 
\begin{equation}\label{lastestimateoftheday}
    \begin{aligned}
         S_n \leq \sum_i \frac{\epsilon}{N} \bigg[ S_{n+1} +& \bigg( \frac{t}{\rho^\lambda}\bigg)^{\frac{1}{2-p}} \bigg] + \sum_i c(\epsilon,\gamma) b^n \bigg( \frac{t}{\rho^{\lambda}}  \bigg)^{\frac{1}{2-p}} + \I \leq \epsilon S_{n+1} + \gamma b^n \bigg{\{} \I+ \bigg( \frac{t}{\rho^{\lambda}} \bigg)^{\frac{1}{2-p}}  \bigg{\}},
    \end{aligned}
\end{equation} and the conclusion follows from the classical iteration of Lemma \ref{iteration}.

\subsection*{Standard Anisotropic Geometry: Proof of Theorem \ref{el1l1}}

\noindent Let $0<t<T$ and $ \rho>0 $ such that the following inclusion is satisfied,
\[\q:=\k_{\rho} \times [0,t] \subset \Omega_T.\]
\noindent To consider intermediate cylinders, for a fixed $\sigma \in (0,1]$ we define  
\[ \q_{\sigma} = \k_{\sigma\rho} \times [0,t] = \prod_i \bigg{\{} |x_i| < (\sigma \rho)^{\frac{p}{p_i}} \bigg{\}} \times [0,t], \qquad \text{and} \qquad \q=\q_1.\]
\noindent Moreover, for such fixed $t,\rho$, we define the quantity
\begin{equation} \label{nusigma}
    \nu_{\Sigma}= \sum_k \bigg( \frac{t}{\rho^p}\bigg)^{\frac{1}{2-p_k}}.
\end{equation}

\begin{lemma}\label{auxiliar1}
Let $u$ be a non-negative local weak super-solution to \eqref{EQ} in $\Omega_T$ and $\sigma \in (0,1)$ a number. Then, there exists a positive constant $\gamma$, depending on the data, such that, either there exists an $ i \in \{1,\dots,N\}$ for which \eqref{ealternative} is valid, or for all $i\in \{1, \dots,N\}$ we have
\begin{equation}\label{firstineq}
\begin{aligned}
    \sum_i \rho^{-\frac{p}{p_i}} \iint_{\q_{\sigma}} & |\partial_i u |^{p_i} \, dxd\tau \leq \frac{\gamma} {(1-\sigma)^p} \sum_i  \bigg(\frac{t}{\rho^{\lambda_i}}\bigg)^{\frac{1}{p_i}} \bigg( S+ \nu_{\Sigma} \rho^N \bigg)^{\frac{2(p_i-1)}{p_i}}, 
    \end{aligned}
\end{equation}
with $\lambda_i=N(p_i-2)+p$ and being 
\[\displaystyle{S= \sup_{0\leq \tau\leq t} \int_{\k_\rho} } u(x,\tau)\, dx.\]
\end{lemma}

\begin{proof}
For $\sigma \in (0,1]$  we consider the cylinders 
\[ \q_{\sigma} = \k_{\sigma\rho} \times [0,t] = \prod_i \left\{ |x_i| < (\sigma \rho)^{\frac{p}{p_i}} \right\} \times [0,t], \qquad \text{and} \qquad \q=\q_1.\]
We use the estimates \eqref{Negative-EE} by testing the equation with \[\varphi_i=\tau^{\frac{1}{p_i}}(u+\nu)^{1-\frac{2}{p_i}} \zeta,\] where $\zeta$ is a cut-off function of the type \eqref{zeta}, defined between $\k_{\sigma \rho}$ and $ \k_{\rho}$, therefore verifying 
\[ \|\partial_i \zeta_i\|_{\infty,\k_{\rho}} \leq \gamma/ [(1-\sigma)\rho]^{\frac{p}{p_i}}  \, .\]
This gives, for all $i \in \{1, \dots, N\}$, the inequalities
\begin{equation} \label{gianni}
    \begin{aligned}
        \iint_{\q_{\sigma}} & |\partial_i u |^{p_i} \tau^{\frac{1}{p_i}} (u+\nu_{\Sigma})^{-\frac{2}{p_i}} \, dxd\tau\\
        &\leq \iint_{\q_{\sigma}}  \bigg( \sum_j |\partial_j u |^{p_j}\bigg) \tau^{\frac{1}{p_i}} (u+\nu_{\Sigma})^{-\frac{2}{p_i}} \, dxd\tau \leq \gamma t^{\frac{1}{p_i}} \int_{\k\times\{t\}} (u+\nu_{\Sigma})^{\frac{2(p_i-1)}{p_i}}\, dx\\
        & \qquad \qquad +     \frac{\gamma}{(1-\sigma)^p \rho^p} \sum_j\bigg[1+ C^{p_j} \rho^p
        \bigg]\iint_{\q}  (u+\nu_{\Sigma})^{p_j-\frac{2}{p_i}} \tau^{\frac{1}{p_i}}\, dxd\tau \\
        &\qquad \qquad + \gamma  \bigg( \sum_j C^{p_j}\bigg)  \iint_{\q} (u+\nu_{\Sigma})^{-\frac{2}{p_i}} \tau^{\frac{1}{p_i}}\, dxd\tau.
    \end{aligned} \end{equation} 
    \noindent 
We estimate the various terms. The first integral on the right-hand side of \eqref{gianni} is manipulated as in \eqref{bacalao} to get
\begin{equation*}
    \begin{aligned}
        t^{\frac{1}{p_i}} \int_{\k} (u+\nu_{\Sigma})^{\frac{2(p_i-1)}{p_i}}\, dx &\leq t^{\frac{1}{p_i}} |\k|^{\frac{2-p_i}{p_i}} \bigg( \sup_{0\leq \tau \leq t}\int_{\k} u(x,\tau)\, dx + \nu_{\Sigma} |\k|   \bigg)^{\frac{2(p_i-1)}{p_i}}\\
        & \leq \gamma t^{\frac{1}{p_i}} \rho^{N(\frac{2-p_i}{p_i})} \bigg( S + \nu_{\Sigma} \rho^N   \bigg)^{\frac{2(p_i-1)}{p_i}}\\
        &  = \gamma \rho^{\frac{p}{p_i}} \bigg( \frac{t}{\rho^{\lambda_i}} \bigg)^{\frac{1}{p_i}} (S+\nu_{\Sigma} \rho^N)^{\frac{2(p_i-1)}{p_i}}.
    \end{aligned}
\end{equation*}
\noindent The second term can be estimated by using that $(u+\nu_{\Sigma})^{p_j-2}<\nu_{\Sigma}^{p_j-2}$ to get for all $i=1, \dots, N$ the inequalities
\begin{equation*}
    \begin{aligned}
            \sum_j \bigg( \frac{[1+C^{p_j}\rho^p]}{\nu_{\Sigma}^{2-p_j}} \bigg) &\iint_{\q} (u+\nu_{\Sigma})^{\frac{2(p_i-1)}{p_i}} \tau^{\frac{1}{p_i}} \, dx d\tau\\
            &\leq \sum_j \bigg( \frac{[1+C^{p_j}\rho^p]}{ \nu_{\Sigma}^{2-p_j}} \bigg) t^{1+\frac{1}{p_i}} \bigg( \sup_{0\leq \tau \leq t} \int_{\k} u(x,\tau) \,dx + \nu_{\Sigma} \rho^N\bigg)^{\frac{2(p_i-1)}{p_i}}\rho^{N(\frac{2-p_i}{p_i})}\\
            &\leq \sum_j \bigg( \frac{[1+C^{p_j}\rho^p]}{\nu_{\Sigma}^{2-p_j}} \bigg) \rho^{\frac{p}{p_i}} \bigg(\frac{t}{\rho^{\lambda_i}} \bigg)^{\frac{1}{p_i}} t \bigg( S+\nu_{\Sigma} \rho^N\bigg)^{\frac{2(p_i-1)}{p_i}}.
    \end{aligned} \end{equation*} \noindent
    Finally the third term on the right-hand side of \eqref{gianni} is estimated, for any $i,j \in \{1,\dots, N\}$, as
\begin{equation*}
    \begin{aligned}
        C^{p_j} \iint_{\q} & \tau^{\frac{1}{p_i}} (u+\nu_{\Sigma})^{-\frac{2}{p_i}} \, dxd\tau\\
        &\leq C^{p_j} \rho^{\frac{p}{p_i}} \bigg( \frac{\rho^{\frac{p}{p_j}}}{\nu_{\Sigma}}\bigg)^{p_j} \bigg(  \frac{t}{\rho^p} \bigg) \nu^{p_j-2} \bigg( \frac{t}{\rho^{\lambda_i}} \bigg)^{\frac{1}{p_i}} (S+\nu_{\Sigma} \rho^N)^{\frac{2(p_i-1)}{p_i}}\\
        &\leq \bigg( \frac{C \rho^{\frac{p}{p_j}}}{\nu_{\Sigma}} \bigg)^{p_j} \rho^{\frac{p}{p_i}} \bigg( \frac{t}{\rho^p} \bigg) \nu_{\Sigma}^{p_j-2} \bigg( \frac{t}{\rho^{\lambda_i}}\bigg)^{\frac{1}{p_i}} (S+ \nu_{\Sigma} \rho^N)^{\frac{2(p_i-1)}{p_i}}.
    \end{aligned}
\end{equation*}

\noindent Collecting everything together we obtain 
\begin{equation}\label{giannino}
    \begin{aligned}
         \iint_{\q_{\sigma}} &  |\partial_i u |^{p_i}\tau^{\frac{1}{p_i}} (u+\nu_{\Sigma})^{-\frac{2}{p_i}} \, dxd\tau\\
         &\leq \gamma \rho^{\frac{p}{p_i}} \bigg( \frac{t}{\rho^{\lambda_i}}\bigg)^{\frac{1}{p_i}}  (S+\nu_{\Sigma} \rho^N)^{\frac{2(p_i-1)}{p_i}} \times \\
         & \qquad \qquad \qquad \qquad \times \bigg{\{} 1+ \sum_j \frac{t}{\nu_{\Sigma}^{2-p_j}\rho^p} [1+C^{p_j} \rho^p] + \sum_j \bigg( \frac{C \rho^{\frac{p}{p_j}}}{\nu_{\Sigma}}\bigg)^{p_j} \bigg( \frac{t}{\rho^p}\bigg) \nu_{\Sigma}^{p_j-2}  \bigg{\}}
    \end{aligned}
\end{equation}
\noindent The second factor on the right of \eqref{giannino} is smaller than $4$ if \eqref{ealternative} is violated for all indexes $j \in \{1, \dots, N\}$, and once we observe 
\[\nu_{\Sigma}= \sum_k \bigg( \frac{t}{\rho^p} \bigg)^{\frac{1}{2-p_k}} \ge \bigg( \frac{t}{\rho^p} \bigg)^{\frac{1}{2-p_j}} \ , \qquad \qquad \forall j=1,\dots,N.\]
\noindent This allows us to evaluate
\begin{equation*}
\begin{aligned}
&  \rho^{-\frac{p}{p_i}} \iint_{\Q_{\sigma \rho}} |\partial_i u |^{p_i-1}\, dx d\tau \\
& = \rho^{-\frac{p}{p_i}} \iint_{\Q_{\sigma \rho}} \bigg( |\partial_i u |^{p_i-1} \tau^{\frac{1}{p_i} (\frac{p_i-1}{p_i})} (u+ \nu)^{\frac{-2}{p_i}(\frac{p_i-1}{p_i})} \bigg) \, \bigg( \tau^{\frac{-1}{p_i} (\frac{p_i-1}{p_i})} (u+ \nu)^{\frac{+2}{p_i}(\frac{p_i-1}{p_i})} \bigg) \, dx d\tau \\
&\leq   \rho^{-\frac{p}{p_i}}\bigg(\iint_{\Q_{\sigma \rho}} |\partial_i u |^{p_i} \tau^{\frac{1}{p_i}} (u+\nu)^{\frac{-2}{p_i}}\, dxd\tau \bigg)^{\frac{p_i-1}{p_i}} \bigg( \iint_{\Q_{\sigma \rho}} \tau^{\frac{-1}{p_i}(p_i-1)} (u+\nu)^{\frac{2}{p_i}(p_i-1)} \, dx d\tau \bigg)^{\frac{1}{p_i}}  \\
& \leq \rho^{-\frac{p}{p_i}}\bigg(  \frac{\gamma}{(1-\sigma)^p} \ \rho^{\frac{p}{p_i}} \bigg( \frac{t}{\rho^{\lambda_i}}\bigg)^{\frac{1}{p_i}}  (S+\nu_{\Sigma} \rho^N)^{\frac{2(p_i-1)}{p_i}}\bigg)^{\frac{p_i-1}{p_i}}\bigg( \gamma \, \bigg( \frac{t}{\rho^{\lambda_i}} \bigg)^{\frac{1}{p_i}} \, \rho^{\frac{p}{p_i}} \bigg{\{} S+\bigg( \frac{t}{\rho^{\lambda}}  \bigg)^{\frac{1}{2-p}}\bigg{\}}^{\frac{2(p_i-1)}{p_i}} \bigg)^{\frac{1}{p_i}}\\
& \leq  \frac{\gamma}{(1-\sigma)^p} \bigg( \frac{t}{\rho^{\lambda_i}}\bigg)^{\frac{1}{p_i}}  (S+\nu_{\Sigma} \rho^N)^{\frac{2(p_i-1)}{p_i}}.
\end{aligned}
\end{equation*}
\end{proof}

\subsection*{Proof of Theorem \ref{el1l1} concluded.}

\begin{proof}
  We fix \(\rho>0\), define the sequence of increasing radii
  \[\ \rho_{n} := \rho \sum_{k=0}^n 2^{-k} , \qquad 
  \rho =\rho_o \leq\rho_n \leq \tilde{\rho_n}:= \frac{\rho_n+\rho_{n+1}}{2} \leq \rho_{n+1} < \rho_\infty= 2\rho\]
and construct the family of concentric standard anisotropic cubes
  \[
  \k_n=  \prod_i \bigg{\{}|x_i|<  \rho_n^{\frac{p}{p_i}}\bigg{\}}, \qquad \tilde{\k}_n=  \prod_i\bigg{\{}|x_i|< \tilde{\rho_n}^{\frac{p}{p_i}}\bigg{\}},\]
  verifying \(\k_n \subset \tilde{\k}_n \subset \k_{n+1}\), and for any $\tau_1, \tau_2 \in [0,t]$, we consider the family of cylinders
  \[\q_n=  \k_n \times [\tau_1, \tau_2] \subset \tilde{\q}_n=  \tilde{\k}_n \times [\tau_1, \tau_2] \subset \q_{n+1}.\]
  \noindent For each $n \in \N\cup\{0\}$ chosen, consider $\zeta_n(x)$ a cut-off function of the form \eqref{zeta} between $\k_n$ and $\tilde{\k}_n$ that is time-independent and verifies 
  \[0 \leq \zeta_n \leq 1, \qquad (\zeta_n)_{| \partial \tilde{\k}_n}=0, \qquad \qquad \|\partial_i \zeta_{n}\|_{\infty,\tilde{\k}_n} \leq \gamma \left(\dfrac{2^{n}}{\rho} \right)^{\frac{p}{p_i}}.\]
  
  \noindent Testing \eqref{EQ}-\eqref{structure-conditions} with such a $\zeta_n$ we obtain 
\begin{equation}\label{start1}
    \begin{aligned}
        \int_{\k_n}& u(x,\tau_1)\, dx\leq \int_{\tilde{\k}_n} u(x,\tau_2)\, dx +  \sum_i  \bigg( \|\partial_i \zeta_n \|_{\infty}  C_1+ C \bigg) \iint_{\tilde{\q}_n}  |\partial_i u |^{p_i-1}\, dxd\tau  \\
        &\qquad \qquad \qquad \qquad + \sum_i (C^{p_i-1}\|\partial_i \zeta_n \|_{\infty} + C^{p_i} ) \iint_{\tilde{\q}_n} \, dxd\tau    \,.
    \end{aligned}
\end{equation} for arbitrary time levels $\tau_1,\tau_2 \in [0,t]$. Again, by the continuity of $u$ as a map $[0,T]\rightarrow L^{2}_{loc}(\Omega)$, we take $\tau_2$ as the time level in $[0,t]$ such that 
\[ \I= \inf_{0\leq \tau \leq t} \int_{\k_{2\rho}} u(x, \tau)\, dx = \int_{\k_{2\rho}} u(x, \tau_2)\, dx,\]
and set 
\[
S_n:= \sup_{0\leq \tau \leq t} \int_{\k_n} u(x,\tau)\, dx \ . \]
Since $\tau_1$ is arbitrary, \eqref{start1} yields
\begin{equation*} \label{cicciolotto} S_{n} \leq \I+ \gamma 2^{\frac{p}{p_1} n }\sum_i \rho^{-\frac{p}{p_i}}\iint_{\tilde{\q}_n} | \partial_i u|^{p_i-1}\,dx\, d\tau + \gamma 2^{\frac{p}{p_1} n } \sum_i \bigg( C^{p_i-1} \rho^{-\frac{p}{p_i}}+C^{p_i} \bigg) \iint_{\tilde{\q}_n} \, dxd\tau .\end{equation*} 
\noindent The last term on the right-hand is dominated as follows:
\begin{equation*}
    \begin{aligned}
        \bigg( C^{p_i-1} \rho^{-\frac{p}{p_i}}+C^{p_i} \bigg) \iint_{\tilde{\q}_n} \, dxd\tau
        &\leq \gamma \bigg[ \bigg( \frac{C \rho^{\frac{p}{p_i}}}{\nu_{\Sigma}}  \bigg)^{p_i-1}+ \bigg(   \frac{C^{p_i \rho^p}}{\nu_{\Sigma}^{p_i-1}}\bigg) \bigg] \bigg(  \sum_j \bigg( \frac{t}{\rho^{\lambda_j}} \bigg)^{\frac{1}{2-p_j}} \bigg)\\
        &\leq \gamma  \sum_j \bigg( \frac{t}{\rho^{\lambda_j}} \bigg)^{\frac{1}{2-p_j}},
    \end{aligned}
\end{equation*} \noindent recalling  $t < \nu_{\Sigma}^{2-p_i} \rho^{p}$, for all $i=1, \dots, N$, and assuming that  condition \eqref{ealternative} is violated for all indexes. Therefore, by applying first Lemma \ref{auxiliar1} to the pair of cylinders $\q_n$ and $\tilde{\q}_n$, for which $1-\sigma \geq 2^{-(n+4)}$, and then Young's inequality one gets
\begin{equation*} 
\begin{aligned}
S_{n} & \leq \I+ \gamma 2^{n \frac{p}{p_1}}\sum_i \rho^{-\frac{p}{p_i}}\iint_{\tilde{\q}_n} | \partial_i u|^{p_i-1}\,dx\, d\tau + \gamma 2^{n \frac{p}{p_1}} \sum_i \bigg( \frac{t}{\rho^{\lambda_i}} \bigg)^{\frac{1}{2-p_i}} \\
& \leq  \I + \gamma b^n \sum_i \left(\frac{t}{\rho^{\lambda_i}} \right)^{\frac{1}{p_i}} \bigg{\{} S_{n+1}+ \nu_{\Sigma} \rho^N  \bigg{\}}^{\frac{2(p_i-1)}{p_i}} + \gamma2^{n\frac{p}{p_1}} \sum_i \bigg( \frac{t}{\rho^{\lambda_i}} \bigg)^{\frac{1}{2-p_i}} \\
 &  \leq \epsilon S_{n+1} + \gamma (\epsilon) b^n \bigg{\{} \I+ \sum_i \bigg( \frac{t}{\rho^{\lambda_i}} \bigg)^{\frac{1}{2-p_i}} \bigg{\}}, \qquad \qquad b>1.
\end{aligned}
\end{equation*} 

\noindent A standard iteration finishes the proof as in the case of \eqref{lastestimateoftheday}
\end{proof}



\section{Proof of the backward $L^r$-$L^{\infty}$ estimates} \label{lrlinfty}

\noindent  The proof of Theorems \ref{lr-decay}-\ref{elr-decay} rely on two estimates: $L^r$-$L^\infty$ estimates combined with a $L^r$ estimates backward in time; the presentation is done separately for the intrinsic and the standard geometries.

     \subsection*{Intrinsic Anisotropic Geometry: Proof of Theorem \ref{lr-decay}}

\noindent

\begin{theorem}[$L^r_{loc}$-$L^{\infty}_{loc}$ estimates]\label{sup}
   Suppose $u$ is a non-negative, locally bounded, local weak sub(super)-solution to \eqref{EQ}-\eqref{structure-conditions} in $\Omega_T$. Let $r \ge 1$ and $\lambda_r= N(p-2)+rp>0$. Then, there exists a positive constant $\gamma$, depending only on the data, such that 
\[\forall t>0, \quad \forall \rho>0 \ : \,\,   \K_{4\rho}(t)\times (0,t) \subset \Omega_T, \]
either \eqref{alternative} holds for some $i \in \{1,\dots, N\}$ or 
\begin{equation}\label{sup-estimate}
    \sup_{\K_{{\rho}/{2}}(t)\times [t/2,t]} u \leq \gamma \bigg(\frac{t}{\rho^p} \bigg)^{\frac{-N}{\lambda_r}} \bigg(\sup_{0\leq \tau \leq t} \dashint_{\K_{\rho}(t)} u^r(x,\tau)\, dx  \bigg)^{\frac{p}{\lambda_r}}+ \gamma \bigg( \frac{t}{\rho^p}\bigg)^{\frac{1}{2-p}}.
\end{equation}

\end{theorem}

\begin{proof}
Assume condition \eqref{alternative} does not hold for every $i \in \{1,\dots, N\}$. Let $\sigma \in (0,1)$ be fixed and consider the decreasing sequences of radii, for each $i \in \{ 1, \dots, N\}$,
\begin{equation*} \label{rhoi}  \rho_i:=\rho^{\frac{p}{p_i}} \bigg(  \frac{t}{\rho^p}\bigg)^{\frac{(p_i-p)}{(2-p)p_i}}, \qquad \qquad  \rho_{n,i}:= \rho_i \left( \sigma + \frac{1-\sigma}{2^n} \right)^{\frac{p}{p_i}}, \end{equation*} and of time levels
\[ \sigma t = t_{\infty} < t_n := t \left( \sigma + \frac{1-\sigma}{2^n} \right) \leq t_0=t \]
from which one constructs the sequence of nested and shrinking cylinders 
\[\Q_n= \K_n \times (t-t_n, t), \qquad \text{for} \qquad \K_n = \prod_i \bigg{\{} |x_i|< \rho_{n,i} \bigg{\}} \,  .\]
\noindent For each $n \in \N$, let $\displaystyle{\zeta_n(x,t) = \prod_i \zeta_i^{p_i}(x_i) \eta(t)}$ be a cut-off function as in \eqref{zeta} therefore verifying
\[\zeta_i(x_i)= \left\{
\begin{array}{cc}
1, & \ |x_i| < \rho_{(n+1),i}\\[.4em]
0, & \ |x_i| \geq  \rho_{n,i}
\end{array} \right. \ , 
\qquad 
\|\partial_i \zeta_i\|_{\infty} \leq \left(\frac{2^{n+1}}{(1-\sigma)\rho}\right)^{\frac{p}{p_i}} \bigg( \frac{t}{\rho^p} \bigg)^{\frac{(p-p_i)}{(2-p)p_i}}, \]
for all $i=1, \dots, N$, and 
\[
\eta(\tau) = \left\{
\begin{array}{cc}
0 & , \  0 \leq \tau \leq t-t_n \\[.4em]
1 & , \  t-t_{n+1}\leq \tau \leq t 
\end{array} \right. \ , 
\qquad |\partial_{\tau}\eta| \leq \frac{2^{n+1}}{(1-\sigma)t} \ .
\]
In the weak formulation \eqref{EQ-weak}, for each $n \in \N$, consider the test function $\varphi_n= (u-k_{n+1})_+ \xi_n$, over the cylinders $\Q_n$, for the truncation levels
\[ 0 \leq k_n= k \left( 1-\frac{1}{2^n}\right) < k \ ,\qquad \qquad  n\in \N \cup \{0\}, \]
where $k$ is a positive real number to be determined. By the classical energy estimate \eqref{general-EE} we obtain the following bound on the energy

\begin{equation*} \label{standard-EE}
    \begin{aligned}
        \E_n:=\sup_{t-t_n\leq \tau \leq t}& \int_{\K_n \times \{\tau\}} (u-k_{n+1})_+^2 \xi_n\, dx+ \sum_i \iint_{\Q_n} |\partial_i [(u-k_{n+1})_+ \xi_{n}]|^{p_i} dx \, d\tau  \\
        &\leq \gamma \|\partial_{\tau} \eta\|_{\infty } \iint_{\Q_n} (u-k_{n+1})_+^2 \, dxd\tau +\\
        & +\gamma \bigg{\{} \sum_i \bigg( \|\partial_i \zeta_i \|_{\infty}^{p_i} + C^{p_i}\bigg)  \iint_{\Q_n} (u-k_{n+1})_+^{p_i} \,dx d\tau + C^{p_i} \iint_{\Q_n} \chi_{[u>k_{n+1}]} \,  dxd\tau\, \bigg{\}}\\
        & \leq \frac{\gamma2^n}{(1-\sigma)t}  \iint_{\Q_n} (u-k_{n+1})_+^2 \, dxd\tau +\\
        & +\gamma  \sum_i \bigg( \frac{2^{np_i}}{(1-\sigma)^p \rho^p} \bigg( \frac{t}{\rho^p}\bigg)^{\frac{p_i-p}{2-p}} + C^{p_i}  \bigg)  \iint_{\Q_n} (u-k_{n+1})_+^{p_i} \,dx d\tau \\
        & + \gamma \sum_i \frac{1}{t} (tC^{p_i}) \iint_{\Q_n} \chi_{[u>k_{n+1}]} \,  dxd\tau\, \\
        & \leq \frac{\gamma 2^{2n}}{(1-\sigma)^pt}  \bigg{\{} \iint_{\Q_n} (u-k_{n+1})_+^2 \, dxd\tau + \sum_i \bigg(\frac{t}{\rho^p} \bigg)^{\frac{2-p_i}{2-p}} \iint_{\Q_n} (u-k_{n+1})_+^{p_i} \,dx d\tau\\
        & + \bigg( \frac{t}{\rho^p} \bigg)^{\frac{2}{2-p}} \iint_{\Q_n} \chi_{[u>k_{n+1}]} \,  dxd\tau\, \bigg{\}},
    \end{aligned}
\end{equation*}\noindent 
where first we implemented the construction of the cut-off function $\zeta$ and then we have used that for each $i\in \{1, \dots, N\}$ the condition \eqref{alternative} is violated.
\noindent \subsection*{The case $\max\left\{1, \frac{2N}{N+2} \right\} < p < 2 $ }
\noindent We estimate the energy $\E_n$ from above in terms of the $L^2$-norm of the truncations $(u-k_n)_+$. Observe that for all $s=0,1, \dots, N$, having defined $p_0=2$, it holds
\begin{equation*} \label{stimette} 
\begin{aligned}
\iint_{\Q_n} (u-k_{n})_+^2\, dxd\tau & \geq   \iint_{\Q_n \cap [u>k_{n+1}]} (u-k_{n})_+^{2-p_s} (u-k_{n})_+^{p_i} \, dxd\tau\\
& \geq \left( \frac{k}{2^{n+1}} \right)^{2-p_s} \iint_{\Q_n \cap [u>k_{n+1}]}  (u-k_{n})_+^{p_s} \, dxd\tau\\
& \geq  \left( \frac{k}{2^{n+1}} \right)^{2-p_s} \iint_{\Q_n }  (u-k_{n+1})_+^{p_s}\, dxd\tau.
\end{aligned}
\end{equation*}
 
\noindent Hence we have 
\begin{equation} \label{cicciotto}
\begin{aligned}
\E_n \leq \frac{\gamma2^{2n}}{(1-\sigma)^pt}  \bigg{\{} 1  + \sum_i \bigg(\frac{t}{\rho^p} \bigg)^{\frac{2-p_i}{(2-p)}} \frac{2^{n(2-p_i)}}{k^{2-p_i}}+ \bigg( \frac{t}{\rho^p} \bigg)^{\frac{2}{2-p}} \frac{2^{2n}}{k^2} \bigg{\}}\iint_{\Q_n} (u-k_n)_+^2 \, dxd\tau,
        \end{aligned} \end{equation}
\noindent and taking into account as a further condition 
    \begin{equation}\label{kappa}
    k\geq \left(\frac{t}{ \rho^p}\right)^{\frac{1}{2-p}} \ , 
    \end{equation}
    \noindent the right hand side of \eqref{cicciotto} now reads
\begin{equation} \label{ricicciotto}
\E_n\leq \frac{\gamma 2^{4n}}{(1-\sigma)^p t} \iint_{\Q_n} (u-k_{n})_+^2 \, dx d\tau.\end{equation}
\noindent Now we want to put in a chain the estimate of $\E_n$ obtained in terms of $\|(u-k_n)_+\|_{L^2(\Q_n)}^2$ with the 
 anisotropic Sobolev embedding \eqref{embedding}.\newline
 \noindent Here we take advantage of exponent $p$ being in the super-critical range, $p> \max\{1,\, 2N/(N+2)\}$: indeed, in such a range, the number $q= p(N+2)/N$ is greater than $2$ and we can use H\"older inequality on $\|(u-k_{n+1})_+\|_{L^2(\Q_{n+1})}^2$ to allow the aforementioned chaining procedure. In the embedding \eqref{embedding} we make the choices 
 \[q=\frac{p(N+2)}{N}, \qquad \text{and} \qquad \theta= \frac{p}{p^*}, \qquad \sigma= 2,\]
 to get
    \begin{equation*}
    \begin{aligned}
      &\iint_{\Q_{n+1}} (u-k_{n+1})_+^2\xi_n^2\,  dxd\tau \\
      & \leq  \bigg( \iint_{\Q_{n}} ((u-k_{n+1})_+\xi_n)^{p(\frac{N+2}{N})} \, dxd\tau\bigg)^{\frac{2N}{p(N+2)}} | \Q_{n} \cap [u>k_{n+1}]|^{1-\frac{2N}{p(N+2)}} \\
      & \leq\gamma \left[ \left(\sup_{t-t_n \leq \tau \leq t} \int_{\K_{n} \times\{\tau\}} (u-k_{n+1})_+^2 \xi^2 \, dx\right)^{p/N} \bigg( \prod_i \iint_{\Q_{n}}| \partial_i((u-k_{n+1})_+ \xi )|^{p_i} \, dxd\tau\bigg)^{\frac{p}{Np_i}}\right]^{\frac{2N}{p(N+2)}} \\[.2em]
      & \qquad \qquad \qquad \qquad \times  | \Q_{n} \cap [u>k_{n+1}]|^{1-\frac{2N}{p(N+2)}} \\[.4em]
      & \leq  \gamma  \bigg[ \E_n^{\frac{p}{N}} \prod_i\E_n^{\frac{p}{Np_i}} \bigg]^{\frac{2N}{p(N+2)}}   | \Q_n \cap [u>k_{n+1}]|^{1-\frac{2N}{p(N+2)}}\\[.4em]
      & \leq  \gamma \E_n^{(\frac{p+N}{N+2})(\frac{2}{p})} \bigg(\frac{2^{2n}}{k^2}\iint_{\Q_n} (u-k_n)_+^2\, dxd\tau  \bigg)^{\frac{N(p-2)+2p}{p(N+2)}}\\[.4em]
      & \leq \frac{\gamma b^n}{[(1-\sigma)^pt]^{{(\frac{N+p}{N+2})(\frac{2}{p})}} k^{(\frac{2}{p})\frac{N(p-2)+2p}{N+2}}} \bigg(\iint_{\Q_n} (u-k_n)_+^2\, dxd\tau  \bigg)^{1+\frac{2}{N+2}},\qquad \text{for} \quad b>1.
    \end{aligned}
    \end{equation*}
\noindent By setting $X_n= |\Q_n|^{-1} \|(u-k_n)_+\|_{2, \Q_n}^2$, from the previous estimate we derive
\begin{equation}\label{Xn}
X_{n+1}\leq \frac{\gamma b^n}{[(1-\sigma)^p]^{{(\frac{N+p}{N+2})(\frac{2}{p})}} k^{(\frac{2}{p})\frac{\lambda_2}{N+2}}} \bigg(\frac{\rho^p}{t}  \bigg)^{\frac{2N}{p(N+2)}} \ X_{n}^{1+\frac{2}{N+2}} ,
\end{equation} \noindent 
with $\lambda_2=N(p-2)+2p$. By choosing $k>0$ such that 
\begin{equation*}\label{secondkest}
\dashiint_{\Q_0} u^2 \leq \gamma^{-\frac{N+2}{2}} b^{-\left(\frac{N+2}{2}\right)^2} (1-\sigma)^{(N+p)} \bigg( \frac{t}{\rho^p}\bigg)^{\frac{N}{p}} k^{\frac{\lambda_2}{p}},
\end{equation*}
the Fast Converge Lemma \ref{fastgeomconv}, ensures $X_n \rightarrow 0 $ as $n \rightarrow \infty$, meaning that
\[ \sup_{\K_{\sigma \rho}(t) \times [\sigma t,\,t]} u \leq k \leq \frac{\gamma}{(1-\sigma)^{\frac{p(N+p)}{\lambda_2}}} \bigg( \frac{t}{\rho^p} \bigg)^{-\frac{N}{\lambda_2}} \bigg( \dashiint_{\K_{ \rho}(t) \times [0,\,t]} u^2\, dxd\tau\bigg)^{\frac{p}{\lambda_2}}+ \bigg( \frac{t}{\rho^p}\bigg)^{\frac{1}{2-p}},\]
\noindent and then
\begin{equation*} \label{rest}
\begin{aligned}
    & \sup_{\K_{\sigma \rho}(t) \times [\sigma t,\,t]} u   \leq    \frac{\gamma}{(1-\sigma)^{\frac{p(N+p)}{\lambda_2}}}  \bigg(\frac{t}{\rho^p} \bigg)^{-\frac{N}{\lambda_2}} \bigg( \dashiint_{\K_{ \rho}(t) \times [0,\,t]} u^2\, dxd\tau\bigg)^{\frac{p}{\lambda_2}}+ \gamma \bigg( \frac{t}{\rho^p}\bigg)^{\frac{1}{2-p}}\\
    & \qquad \qquad \leq \frac{\gamma}{(1-\sigma)^{\frac{p(N+p)}{\lambda_2}}} \bigg(\frac{t}{\rho^p} \bigg)^{-\frac{N}{\lambda_2}} \bigg(\sup_{\K_{ \rho}(t) \times [0,t]} u \bigg)^{\frac{p(2-r)}{\lambda_2}} \bigg(\dashiint_{\K_{ \rho}(t) \times [0,\,t]} u^r\, dxd\tau\bigg)^{\frac{p}{\lambda_2}}+ \gamma \bigg( \frac{t}{\rho^p}\bigg)^{\frac{1}{2-p}} \ ,
\end{aligned} \end{equation*} 
for every $1 \leq r \leq 2 < q$ for which (and for sure) $\lambda_r=N(p-2)+rp >0$.

\noindent Here we observe that {\it{a priori}} information on the boundedness of $u$ was not necessary in order to get the first sup-estimate in this case. 

\noindent Finally, we perform a cross-iteration on $\sigma \in (0,1)$ as follows. Still referring to radii $\rho_i$ as in the construction above, we now consider the increasing sequences, for $n \in \N \cup \{0 \}$,
\[\tilde{\rho}_{0,i} =\sigma \rho_i \ ,  \qquad \tilde{\rho}_{n,i} =\rho_i \left(\sigma + (1-\sigma) \sum_{j=1}^{n} 2^{-j}\right), \]
\[ \tilde{t}_0 = \sigma t \ ,  \qquad \tilde{t}_n =t \left(\sigma + (1-\sigma) \sum_{j=1}^{n} 2^{-j}\right),\]
\[ \tilde{\K}_n= \prod_i \bigg{\{} |x_i|<  \tilde{\rho}_{n,i}\bigg{\}}, \qquad \tilde{\Q}_n = \tilde{\K}_n \times (t-\tilde{t}_n, t),\]
and define
\[S_n= \sup_{\tilde{\Q}_n} u \ . \]
The previous estimate applied to the pair of cylinders $\tilde{\Q}_n$ and $\tilde{\Q}_{n+1}$ gives us
\begin{eqnarray*}
    S_n & \leq & \frac{\gamma}{(1-\sigma)^{\frac{p(N+p)}{\lambda_2}}} S_{n+1}^{\frac{p(2-r)}{\lambda_2}} \bigg(\frac{t}{\rho^p}\bigg)^{-\frac{N}{\lambda_2}}\bigg( \dashiint _{\tilde{\Q}_{n+1}} u^r\, dxd\tau\bigg)^{\frac{p}{\lambda_2}}+ \gamma \bigg( \frac{t}{\rho^p}\bigg)^{\frac{1}{2-p}}\\
    & \leq & \frac{1}{2} S_{n+1} + \frac{\gamma}{(1-\sigma)^{\frac{p(N+p)}{\lambda_r}}}
    \bigg(\frac{t}{\rho^p}\bigg)^{-\frac{N}{\lambda_r}} \left(\dashiint_{\tilde{\Q}_{\infty}} u^r\, dxd\tau\right)^{\frac{p}{\lambda_r}} + \gamma \left(\frac{t}{\rho^p}\right)^{\frac{1}{2-p}}
\end{eqnarray*}
by means of Young's inequality with $\epsilon=1/2$ for exponents $\mu= \frac{\lambda_2}{p(2-r)}>1$ and $\mu'=\lambda_2/\lambda_r$. Therefore, by iteration, one gets
\[
S_0 \leq \left(\frac{1}{2}\right)^n S_{n}+ \left(\sum_{j=0}^{n-1} 2^{-j}\right) 
\frac{\gamma}{(1-\sigma)^{\frac{p(N+p)}{\lambda_r}}}\bigg(\frac{t}{\rho^p}\bigg)^{-N/\lambda_r}\left(\dashiint_{\tilde{\Q}_{\infty}} u^r\, dxd\tau\right)^{\frac{p}{\lambda_r}} + \gamma \left(\frac{t}{ \rho^p}\right)^{\frac{1}{2-p}} \ . \]
and, by taking $\sigma = 1/2$ and letting $n \rightarrow \infty$ 
\[\sup_{\K_{\frac{\rho}{2}}(t) \times \left[ t/2, t\right]} u \, =\,  \sup_{\tilde{\Q}_{o}} u   \leq  \gamma \bigg(\frac{t}{\rho^p}\bigg)^{-N/\lambda_r} \left(  \dashiint_{\K_{\rho} \times [0,t]} u^r \, dxd\tau \right)^{\frac{p}{\lambda_r}} +  \gamma \left(\frac{t}{\rho^p} \right)^{\frac{1}{2-p}} \ . \]

\subsection*{The case $1 < p \leq \max\left\{1, \frac{2N}{N+2} \right\} $ }
\noindent In this case, the conditions $\lambda_r>0$ and $1<p \leq 2N/(N+2)$ imply $r>2$ and also $q= p\frac{N+2}{N} \leq 2 < r$. Here we need to consider the $L^r$-norm of the truncated functions
\[Y_n= \iint_{\Q_n} (u-k_n)_+^r\, dxd\tau,\]
and supposing $u$ locally bounded, recalling $q<2<r$, we apply the anisotropic embedding \eqref{embedding} to get
\begin{equation*} \begin{aligned}
    &Y_{n+1}  \leq  \iint_{\Q_n} (u-k_{n+1})_+^{r-q} (u-k_{n+1})_+^q \xi^q_n \, dxd\tau  \\
    & \leq  \bigg(\sup_{\Q_0} u\bigg)^{r-q} \iint_{\Q_n} (u-k_{n+1})_+^q \xi^q_n\, dxd\tau\\
    & \leq  \gamma \bigg(\sup_{\Q_0} u\bigg)^{r-q} \bigg( \sup_{t-t_n \leq \tau \leq t}\int_{\K_n} (u-k_{n+1})_+^2 \xi^2_n \, dx \bigg)^{\frac{p}{N}} \bigg(  \prod_i \iint_{\Q_n}| \partial_i((u-k_{n+1})_+ \xi_n |^{p_i}\, dxd\tau\bigg)^{\frac{p}{Np_i}}\\
    &\leq   \gamma \left(\sup_{\Q_0} u\right)^{r-q}  \E_n^{1+\frac{p}{N}}.
\end{aligned} \end{equation*}
\noindent Now again we make a chain of inequalities, but this time using  $\E_n$ and $Y_n$. By acting in a similar fashion as before and assuming \eqref{kappa}, we get
\begin{equation*}
\E_n\leq  \frac{\gamma 2^{n (r+2)}}{(1-\sigma)^p t} \frac{1}{k^{r-2}}\, \, Y_n,\end{equation*} and 
\noindent therefore the aforementioned chain reads
\begin{equation*}
    Y_{n+1} \leq \gamma \bigg( \sup_{\Q_0} u \bigg)^{r-q} \frac{b^n}{((1-\sigma)^p t)^{\frac{(N+p)}{N}} k ^{\frac{(r-2)(N+p)}{N}}} \ Y_n^{1+\frac{p}{N}} \, , \qquad \quad  b= 2^{\frac{(r+2)(N+p)}{N}}>1.
\end{equation*}
\noindent Again by the Fast Convergence Lemma \ref{fastgeomconv}, if $k>0$ is taken so that
\[ Y_0 \leq \gamma^{-\frac{N}{p}} b^{-\frac{N^2}{p^2}} \left(\sup_{\Q_0} u\right)^{-\frac{(r-q)N}{p}} ((1-\sigma)^p t)^{\frac{(N+p)}{p}} k^{\frac{(r-2)(N+p)}{p}}, \]
we obtain $u < k$ for almost every $(x,\tau) \in \Q_{\infty}$. Therefore we choose
\begin{equation} \label{choicek}
k= \gamma \left(\sup_{\Q_0} u\right)^{\frac{(r-q)N}{(N+p)(r-2)}} \left(\iint_{\Q_0 } u^r \, dxd\tau\right)^{\frac{p}{(N+p)(r-2)}} ((1-\sigma)^p t)^{-\frac{1}{r-2}} +  \left(\frac{t}{\rho^p}\right)^{\frac{1}{2-p}}
\end{equation}
for which we get
\[\sup_{\Q_\infty} u \leq \gamma \left(\sup_{\Q_0} u\right)^{\frac{(r-q)N}{(N+p)(r-2)}} \frac{1}{((1-\sigma)^p t)^{\frac{1}{r-2}} }\left(\iint_{\Q_0}  u^r\, dxd\tau \right)^{\frac{p}{(N+p)(r-2)}}  + \left(\frac{t}{\rho^p}\right)^{\frac{1}{2-p}} \ .\]
Proceeding as before, one has
\begin{eqnarray*}
    S_n & \leq & S_{n+1}^{\frac{(r-q)N}{(N+p)(r-2)}} \frac{\gamma}{\left((1-\sigma)^p t\right)^{\frac{1}{r-2}}}\left( \iint_{\tilde{Q}_{n+1}} u^r \, dxd\tau\right)^{\frac{p}{(N+p)(r-2)}}+ \left(\frac{t}{\rho^p}\right)^{\frac{1}{2-p}}\\
    & \leq & \frac{1}{2} S_{n+1} + \frac{\gamma}{\left((1-\sigma)^p t\right)^{\frac{N+p}{\lambda_r}}}\left(\iint_{\tilde{Q}_{\infty}} u^r \, dxd\tau\right)^{\frac{p}{\lambda_r}} + \left(\frac{t}{ \rho^p}\right)^{\frac{1}{2-p}}
\end{eqnarray*}
by means of Young's inequality with $\epsilon=1/2$ for exponent $\mu= \frac{(N+p)(r-2)}{N(r-q)}>1$. Then by iteration, taking $\sigma = 1/2$ and letting $n \rightarrow \infty$
\begin{equation*}
    \begin{aligned}
        \sup_{\K_{\rho/2}(t) \times [t/2, t]} \,  u \, \, &\leq \,    \gamma \ t^{-\frac{N+p}{\lambda_r}} \bigg(  \int_0^t \int_{\K_{\rho}(t) }  u^r  \, dxd\tau\bigg)^{\frac{p}{\lambda_r}} +  \bigg(\frac{t}{\rho^p} \bigg)^{\frac{1}{2-p}}\\
        &=  \gamma  \bigg(\frac{t}{\rho^p}\bigg)^{-\frac{N}{\lambda_r}} \left(  \dashint_0^t \dashint_{\K_{\rho}(t)} u^r \, dxd\tau\right)^{\frac{p}{\lambda_r}} +  \left(\frac{t}{\rho^p} \right)^{\frac{1}{2-p}}\, .\end{aligned}
\end{equation*}

\end{proof}

\begin{theorem}[$L^r_{loc}$ estimates backward in time] \label{Y}
 Let $u$ be a non-negative, locally bounded, local weak solution to \eqref{EQ}-\eqref{structure-conditions} and assume $u \in L^r_{loc}(\Omega_T)$, for some $r>1$. Then there exists a positive constant $\gamma$, depending only on the data, such that either \eqref{alternative} is satisfied for some $i\in \{1, \dots,N\}$ or  
 \begin{equation}\label{lr-linfty}
     \sup_{0\leq \tau \leq t} \int_{\K_{\rho}(t)} u^r(x,\tau)\, dx \leq \gamma \int_{\K_{2\rho}(t)} u^r(x,0)\, dx + \gamma \bigg( \frac{t^r}{\rho^{\lambda_r}} \bigg)^{\frac{1}{2-p}} \, ,
 \end{equation}
\noindent being $\lambda_r=N(p-2)+pr$.
\end{theorem}

\begin{proof}
Assume \eqref{alternative} fails to happen for all $i\in \{1, \dots,N\}$. Fix $\sigma \in (0,1)$ and construct the cylinders
    \[Q_1= \K_{\rho}(t) \times [0,t], \qquad Q_2= \K_{(1+\sigma)\rho}(t)  \times [0,t]. \]
With these stipulations, a cut off function $\zeta$, such as in \eqref{zeta}, between $\K_{\rho}(t)$ and $\K_{(1+\sigma)\rho}(t)$ satisfies
\[ \|\partial_i \zeta_i\|_{\infty} \leq \frac{1}{(\sigma\rho)^{\frac{p}{p_i}}} \bigg( \frac{t}{\rho^p} \bigg)^{\frac{(p-p_i)}{p_i(2-p)}}=:\frac{1}{\sigma^{\frac{p}{p_i}}\rho_i(t)},\]
\noindent and the estimates \eqref{parabolic-estimates} with $K_1= \K_{\rho}(t)$ and $K_2=\K_{(1+\sigma)\rho}(t)$ are now written
\begin{equation*}\label{parabolic-wg} 
\begin{aligned}
\sup_{0\leq \tau \leq t} &\int_{\K_{\rho}(t)} u^r(x, \tau) \, dx \leq \gamma \int_{\K_{(1+\sigma)\rho}(t)} u^r(x,0)\, dx \\
&+ \sum_{i} \frac{\gamma}{\sigma^p\rho^p} \bigg(\frac{t}{\rho^p} \bigg)^{\frac{p-p_i}{2-p}} \bigg{\{} \int_0^t \int_{\K_{(1+\sigma)\rho}(t)} u^{r+p_i-2} \, dxd\tau + 
\\& + \bigg[ \bigg(C \rho^\frac{p}{p_i} \bigg( \frac{t}{\rho^p}\bigg)^{\frac{p_i-p}{p_i(2-p)}}\bigg)^{p_i-1} + \bigg(C \rho^\frac{p}{p_i} \bigg( \frac{t}{\rho^p}\bigg)^{\frac{p_i-p}{p_i(2-p)}}\bigg)^{p_i}  \bigg( 1+ \frac{1}{M_r}\bigg)\bigg] \int_0^t \int_{\K_{(1+\sigma)\rho}(t)} u^{r-1} \, dxd\tau \bigg{\}},
\end{aligned}
\end{equation*} being 
\begin{equation*}\label{Mr-intrinsic}
   M_r= \bigg( \sup_{0\leq \tau \leq t} \dashint_{\K_{\rho}(t)} u^r(x, \tau) \, dx\bigg)^{\frac{1}{r}}.
\end{equation*}
\noindent Without loss of generality one can assume that, for all $i=1, \dots, N$, 
\[C \rho^{\frac{p}{p_i}} \bigg(\frac{t}{\rho^p}\bigg)^{\frac{p_i-p}{p_i(2-p)}}\leq M_r.
\]
\noindent In fact, if for some index $i=1, \dots, N$
\[C \rho^{\frac{p}{p_i}} \bigg(\frac{t}{\rho^p}\bigg)^{\frac{p_i-p}{p_i(2-p)}}> M_r, 
\]
implying that  
\[\sup_{0\leq \tau \leq t} \int_{\K_{\rho}(t)} u^r(x, \tau)\, dx < 
2^N \rho^N \left(C \rho^{\frac{p}{p_i}} \nu^{\frac{p_i-p}{p_i}}\right)^r <
2^N \rho^N  \left( \nu^{\frac{p}{p_i}} \nu^{\frac{p_i-p}{p_i}}\right)^r = \gamma \bigg( \frac{t^r}{\rho^{\lambda_r}} \bigg)^{\frac{1}{2-p}}\]
and then \eqref{lr-linfty} comes immediately. Hence
\begin{equation*}
\begin{aligned}
\sup_{0\leq \tau \leq t}& \int_{\K_{\rho}(t)} u^r(x, \tau) \, dx \leq \gamma \int_{\K_{(1+\sigma)\rho}(t)} u^r(x,0)\, dx +\\
&+ \gamma \sum_{i} \frac{\gamma}{\sigma^p\rho^p} \bigg(\frac{t}{\rho^p} \bigg)^{\frac{p-p_i}{2-p}} \bigg{\{} \int_0^t \int_{\K_{(1+\sigma)\rho}(t)} u^{r+p_i-2} \, dxd\tau + M_r^{p_i-1}\int_0^t \int_{\K_{(1+\sigma)\rho}(t)} u^{r-1} \, dxd\tau \bigg{\}}.
\end{aligned}
\end{equation*}

\noindent We estimate the second integral on the right-hand side  by applying H\"older's inequality,
\begin{equation*} \label{rr}
    \begin{aligned}
         \sum_i \frac{t}{\rho^p} \bigg(\frac{t}{\rho^p} \bigg)^{\frac{(p_i-p)}{(2-p)}} &  \bigg( \sup_{0\leq \tau \leq t} \int_{\K_{(1+\sigma)\rho}(t)} u^{p_i+r-2}(x,\tau)\, dx \bigg)\\
         &\leq \gamma \sum_i \bigg( \frac{t}{\rho^p} \bigg)^{\frac{2-p_i}{2-p}}\bigg( \sup_{0\leq \tau \leq t} \int_{\K_{(1+\sigma)\rho}(t)} u^r(x,\tau)\, dx \bigg)^{\frac{p_i+r-2}{r}} \rho^{\frac{N(2-p_i)}{r}}\\
         &= \gamma \sum_i \bigg( \frac{t^r}{\rho^{\lambda_r}} \bigg)^{\frac{2-p_i}{r(2-p)}}\bigg( \sup_{0\leq \tau \leq t} \int_{\K_{(1+\sigma)\rho}(t)} u^r(x,\tau)\, dx \bigg)^{\frac{p_i+r-2}{r}}.
         \end{aligned}
     \end{equation*} 
     
\noindent The last integral on the right-hand side is dominated as follows
     \begin{equation*}\label{rrr}
         \begin{aligned}
             \gamma &\sum_{i} \frac{\gamma}{\sigma^p\rho^p} \bigg(\frac{t}{\rho^p} \bigg)^{\frac{p-p_i}{2-p}} M_r^{p_i-1}\int_0^t \int_{\K_{(1+\sigma)\rho}(t)} u^{r-1} \, dxd\tau\\
             &  \leq \gamma \sum_{i} \frac{\gamma}{\sigma^p\rho^p} \bigg(\frac{t}{\rho^p} \bigg)^{\frac{p-p_i}{2-p}} M_r^{p_i-1} t \bigg( \sup_{0\leq \tau \leq t} \int_{\K_{(1+\sigma)\rho}(t)} u^{r}(x,\tau) \, dx \bigg)^{\frac{r-1}{r}} (2\rho)^{\frac{N}{r}} \\
             &  \leq \frac{\gamma}{\sigma^p} \sum_i \bigg( \frac{t^r}{\rho^{\lambda_r}} \bigg)^{\frac{2-p_i}{r(2-p)}}\bigg( \sup_{0\leq \tau \leq t} \int_{\K_{(1+\sigma)\rho}(t)} u^r(x,\tau)\, dx \bigg)^{\frac{p_i+r-2}{r}}     \end{aligned}
     \end{equation*} using H\"older inequality and noticing that  
     \[M_r= \bigg( \sup_{0\leq \tau \leq t} \dashint_{\K_{\rho}(t)} u^r(x,\tau)\,  dx \bigg)^{\frac{1}{r}} < \bigg( \sup_{0\leq \tau \leq t} \int_{\K_{\rho(1+\sigma)}(t)} u^r(x,\tau) dx \bigg)^{\frac{1}{r}} (2 \rho)^{-\frac{N}{r}} .\]
    \noindent Putting the estimates all together we finally get
\begin{equation}\label{w} \begin{aligned}
\sup_{0\leq \tau \leq t} &\int_{\K_\rho(t)} u^r (x, \tau) \, dx \, \\
&\leq \gamma\int_{\K_{(1+\sigma)\rho}(t)} u^r(x,0)\, dx +\sum_i  \frac{\gamma}{\sigma^p} \left(\sup_{0\leq \tau \leq t} \int_{\K_{(1+\sigma)\rho}(t)} u^r (x,\tau)\, dx \right)^{\frac{p_i+2-r}{r}} \left(\frac{t^r}{\rho^{\lambda_r}} \right)^{\frac{2-p_i}{r(2-p)}} \ .
\end{aligned} \end{equation}
Now we perform an iteration on $\sigma$: fix $\rho>0$ and for $n\in \N \cup\{0\}$ consider the increasing sequence of radii 
\[\rho_i(t) \leq \rho_{n,i}:= \rho_i(t)\sum_{s=0}^n 2^{-s} \quad \mathrm{so} \ \mathrm{that} \quad \rho_{n+1,i}=(1+\sigma_n) \rho_{n,i}, \ \ \ \mathrm{for} \ \ \sigma_n= \dfrac{\rho_{n+1,i}-\rho_{n,i}}{\rho_{n,i}} \geq \frac{1}{2^{n+2}}. \]
By setting \[S_n = \sup_{0\leq \tau \leq t} \int_{\K_{\rho_n}(t)} u^r(x,\tau) \,dx \ , \]  estimate \eqref{w} now reads
\[S_n \leq \int_{\K_{2\rho}(t)} u^r(x,0) \, dx+\gamma \sum_i  2^{np} \left(S_{n+1}  \right)^{\frac{p_i+r-2}{r}} \left(\frac{t^r}{\rho^{\lambda_r}} \right)^{\frac{2-p_i}{r(2-p)}}.  \]
\noindent We use Young's inequality in each i-th term of the sum 
\begin{equation*}
    \begin{aligned}
         \left[ \gamma 2^{np} \bigg( \frac{t^r}{\rho^{\lambda_r}} \bigg)^{\frac{2-p_i}{r(2-p)}} \right] \bigg( S_{n+1} \bigg)^{\frac{p_i+r-2}{r}}     \leq \epsilon S_{n+1}+  \gamma(\epsilon) b^n \bigg( \frac{t^r}{\rho^{\lambda_r}} \bigg)^{\frac{1}{2-p}},
    \end{aligned}
\end{equation*}
for a constant $b>1$ depending only on the data, and  with these stipulations we arrive at
\[S_n \leq \epsilon \  S_{n+1}  + \gamma(\epsilon) b^n \left(\int_{\K_{2\rho}(t)} u^r(x,0)\, dx + \left(\frac{t^r}{\rho^{\lambda_r}} \right)^{\frac{1}{2-p}}  \right) \, . \]
 \noindent A simple iteration shows
\[S_0 \leq \epsilon^n \  S_{n}  + \gamma(\epsilon) \sum_{k=1}^{n-1} \left( \epsilon \ b \right)^k  \left(\int_{\K_{2\rho}(t)} u^r(x,0)\, dx +\left(\frac{t^r}{\rho^{\lambda_r}} \right)^{\frac{1}{2-p}}  \right),  \]
and proof is completed once we choose $\epsilon= 1/2b <1$ and let $n \rightarrow \infty$ as usual.
\end{proof}
\begin{remark}
Here the exponent $\lambda_r=N(p-2)+pr$ can be of either sign.
\end{remark}
\subsection{Proof of Theorem \ref{lr-decay} concluded}

\begin{proof} We plug inequality \eqref{lr-linfty} into \eqref{sup-estimate} to obtain
\begin{equation*}
        \begin{aligned}
            \|u\|_{\infty, \K_{\rho}(t)\times [t/2,t]} &\leq \gamma t^{-\frac{N}{\lambda_r}} \bigg(   \int_{\K_{2\rho}(t)} u^r(x,0)\, dx + \bigg( \frac{t^r}{\rho^{\lambda_r}} \bigg)^{\frac{1}{2-p}}\bigg)^{\frac{p}{\lambda_r}}+\gamma \bigg( \frac{t}{\rho^p} \bigg)^{\frac{1}{2-p}} \\&
            \leq \gamma  t^{-\frac{N}{\lambda_r}} \bigg( \int_{\K_{2\rho}(t)} u^r(x,0)\, dx\bigg)^{\frac{p}{\lambda_r}}+ \gamma  \bigg( \frac{t}{\rho^p} \bigg)^{\frac{1}{2-p}} .
        \end{aligned}
    \end{equation*}
\end{proof}

     \subsection*{Standard Anisotropic Geometry: Proof of Theorem \ref{elr-decay}}

\noindent

\begin{theorem}{\bf{($L^r_{loc}$- $L^\infty_{loc}$ estimates)}} \label{lrlinfestimate} Let $u$ be a non-negative, locally bounded, local weak sub(super)-solution to \eqref{EQ}-\eqref{structure-conditions} in $\Omega_T$. Let  $r\geq 1$ be such that 
\begin{equation}\label{lambdar}
\lambda_r=N(p-2)+rp >0 \ .
\end{equation} \noindent Then there exists a positive constant $\gamma$, depending only on the data such that, for all $\k_{\rho} \times [0,t] \subset \Omega_T$, either for some $i \in \{1, \dots, N\}$ condition \eqref{ealternative} is satisfied or 
\begin{equation}\label{estLrlinfty}
\sup_{\k_{\rho/2} \times [t/2, t]} u\,   \leq  \gamma \bigg(\frac{t}{\rho^p}\bigg)^{-\frac{N}{\lambda_r}} \bigg(  \sup_{0\leq \tau\leq t }\dashint_{\k_{\rho} } u^r(x,\tau)\, dx \bigg)^{\frac{p}{\lambda_r}} +  \sum_i\left(\frac{t}{\rho^p} \right)^{\frac{1}{2-p_i}} \ .
\end{equation}    
\end{theorem}

\begin{proof}
Assume condition  \eqref{ealternative} is violated for all indexes $i \in \{1, \dots, N\}$. Let $\sigma \in (0,1)$  be fixed and consider the decreasing sequences
\[ \sigma \rho = \rho_{\infty} < \rho_n = \rho \left( \sigma + \frac{1-\sigma}{2^n} \right) \leq \rho_0=\rho\]
and
\[ \sigma t = t_{\infty} < t_n = t \left( \sigma + \frac{1-\sigma}{2^n} \right) \leq t_0=t \]
from which one constructs the sequence of nested and shrinking cylinders 
\[\q_n= \k_n \times (t-t_n, t) \]
where, as usual in the standard anisotropic geometry,
\[\k_n = \prod_i\left\{ |x_i|< \rho_n^{\frac{p}{p_i}} \right\} . \]
Define cutoff function $\zeta_n(x,t)= \zeta_n(x) \xi(\tau)$, as in \eqref{zetatime}, verifying
\[\zeta_{n,i}(x_i)= \left\{
\begin{array}{cc}
1 & , \ |x_i| < \rho_{n+1}\\[.4em]
0 & , \ |x_i| \geq  \rho_{n}
\end{array} \right. \ , 
\qquad 
\|\partial_i \zeta_n \|_{\infty} \,  \leq \left(\frac{2^{n+1}}{(1-\sigma) \rho}\right)^{\frac{p}{p_i}} \]
and 
\[
\xi(\tau)= \left\{
\begin{array}{cc}
0 & , \  0 \leq \tau \leq t-t_n \\[.4em]
1 & , \  t-t_{n+1}\leq \tau \leq t 
\end{array} \right. \ , 
\qquad \|\partial_t \xi \|_{\infty} \leq \frac{2^{n+1}}{(1-\sigma)t} \ .
\]
In the weak formulation \eqref{EQ-weak} we consider test functions  $\varphi_n= (u-k_{n+1})_+ \zeta_n$, over the cylinders $\q_n$, for the truncation levels
\[ 0 \leq k_n= k \left( 1-\frac{1}{2^n}\right) < k \ , n=0,1, \cdots\]
where $k$ is a positive real number to be determined (along the proof). By the energy estimates \eqref{general-EE} we get 

\begin{equation}\label{ugo}
    \begin{aligned}
        \E_n= \sup_{t-t_n \leq \tau \leq t}& \int_{\k_n \times\{\tau\}} (u-k_{n+1})_+^2 \zeta_n \, dx + \sum_i \iint_{\q_n}| \partial_i \left((u-k_{n+1})_+ \zeta_n  \right)|^{p_i}\, dx d\tau \\
         & \leq  \gamma    \frac{2^n}{(1-\sigma)t} \iint_{\q_n} (u-k_{n+1})_+^2 \, dxd\tau \\
        & +\gamma \frac{2^{np}}{(1-\sigma)^p\rho^p} \sum_i\bigg( 1+(C^{p_i}\rho^p )\bigg) \iint_{\q_n} (u-k_{n+1})_+^{p_i}\, dxd\tau \\
        & + \gamma \sum_i C^{p_i} \iint_{\q_n} \chi_{[u>k_{n+1}]}\, dxd\tau .
    \end{aligned}
\end{equation} 

\noindent As in the proof of Theorem \ref{sup}, from now on we distinguish between the case where $p$ is in the super and the sub-critical ranges. We will only present how to proceed when $p$ is in the super-critical range; the sub-critical range is treated analogously to what was done for the anisotropic intrinsic geometry but now taking into account take we are working under the assumptions related to the anisotropic standard setting.

\vskip0.2cm

\noindent Consider $\max \{1, \frac{2N}{N+2} \} < p <2$. By observing that $ \rho^p C^{p_i} \leq 1$, for all $i \in \{1, \dots, N\}$,
\[ \iint_{\q_n} (u-k_{n})_+^2\, dx d\tau \geq \left(\frac{k}{2^{n+1}}\right)^2 \iint_{\q_n} \chi_{[u>k_{n+1}]}\, dx d\tau \ , \]
\[ \iint_{\q_n} (u-k_{n})_+^2 \, dx d\tau\geq \left(\frac{k}{2^{n+1}}\right)^{2-p_i} \iint_{\q_n}  (u-k_{n+1})_+^{p_i}\, dx d\tau\ , \]
 and choosing $k \geq \nu_{\Sigma}$, from the previous estimate \eqref{ugo} one gets
\begin{eqnarray*}
\E_n & \leq & \gamma \frac{2^{(p+2)n}}{(1-\sigma)^p t} \left\{ 1 + \frac{t}{\rho^p} \sum_i k^{p_i-2} + \frac{t}{\rho^p} \sum_i \frac{\rho^p C^{p_i}}{k^2} \right\} \iint_{\q_n} (u-k_{n})_+^2\, dx d\tau \\
& \leq & \gamma \frac{2^{(p+2)n}}{(1-\sigma)^p t} \iint_{\q_n} (u-k_{n})_+^2 \, dx d\tau \ . 
 \end{eqnarray*}
Although the geometry is different, we derive a similar estimate to \eqref{Xn} by means of H\"older's inequality, so to obtain
\begin{equation*} \begin{aligned} &\sup_{\K_{\sigma \rho}(t) \times [\sigma t,\,t]} u\, \, \leq    \\     &\qquad \leq \frac{\gamma}{(1-\sigma)^{\frac{p(N+p)}{\lambda_2}}} \bigg(\frac{t}{\rho^p} \bigg)^{-\frac{N}{\lambda_2}} \bigg(\sup_{\K_{ \rho}(t) \times [0,t]} u \bigg)^{\frac{p(2-r)}{\lambda_2}} \bigg(\dashiint_{\K_{ \rho}(t) \times [0,\,t]} u^r\, dxd\tau\bigg)^{\frac{p}{\lambda_2}}+ \gamma \sum_i \bigg( \frac{t}{\rho^p}\bigg)^{\frac{1}{2-p_i}} \end{aligned} \end{equation*}
An analogous iteration procedure is applied considering the radius to be $\rho$ rather than $\rho_i$, completing thereby the proof for the super-critical range of $p$.

\end{proof}

\begin{theorem}
    {\bf{($L^r_{loc}$ estimates backward in time)}} \label{lrestimate} Let $u$ be a non-negative, locally bounded, local weak solution to \eqref{EQ}-\eqref{structure-conditions} in $\Omega_T$. Assume that $u\in L^r_{loc}(\Omega_T) $, for some $r>1$. Then there exists a positive constant $\gamma$, depending on the data, such that either \eqref{ealternative} is verified for some $i \in \{1,\dots, N\}$, or 
\begin{equation}\label{estimatesLr}
\sup_{0\leq \tau \leq t} \int_{\k_\rho} u^r (x, \tau) \, dx \leq  \gamma  \int_{\k_{2\rho}} u^r(x,0) \, dx+ \gamma \sum_{i}\left(\frac{t^r}{\rho^{\lambda_{i,r}}} \right)^{\frac{1}{2-p_i}},
\end{equation}    \noindent where $\lambda_{i,r}=N(p_i-2)+pr.$
\end{theorem}

\begin{proof}
Assume \eqref{ealternative} is not verified for all $i \in \{1,\dots, N\}$. Fix $\sigma \in (0,1)$ and construct the cylinders
    \[\q= \k_\rho \times [0,t]=\prod_{i} \bigg{\{} |x_i|< \rho^{\frac{p}{p_i}} \bigg{\}} \times [0,t], \qquad \q_{\sigma}= \k_{(1+\sigma)\rho} \times [0,t].\]

\noindent 
Using \eqref{parabolic-estimates} with $Q_1=\q$ and $Q_2=\q_{\sigma}$, and a  time-independent cut-off function $\zeta$ is as in \eqref{zeta} defined in $\K_{(1+\sigma)\rho}$ and verifying  
    \[  \|\partial_i \zeta\|_{\infty} \leq \gamma/(\sigma \rho)^{\frac{p}{p_i}}, \quad \text{for all} \quad i=1, \dots, N,\]
while considering 
\begin{equation}\label{Mr-standard}
   M_r= \bigg( \sup_{0\leq \tau \leq t} \dashint_{\k_{\rho}} u^r \, dx\bigg)^{1/r} > C \rho^{\frac{p}{p_i}} \ , \qquad \forall i=1, \cdots,N
\end{equation}
we obtain 
\begin{eqnarray} \label{estsigma}
\sup_{0\leq \tau \leq t} \int_{\k_{\rho}} u^r(x, \tau) \, dx & \leq & \gamma \int_{\k_{\sigma \rho}} u^r(x,0)\, dx \\ \nonumber
& & + \frac{\gamma}{(\sigma \rho)^p} \sum_{i}  \bigg(1+(C^{p_i} \rho^p) \bigg)  \iint_{\q_{\sigma}} u^{r+p_i-2} \, dxd\tau  \\ \nonumber
& & +\frac{\gamma}{(\sigma \rho)^p}  \sum_{i}  \bigg[(C \rho^{\frac{p}{p_i}})^{p_i-1}  +C^{p_i}\rho^p \bigg(1+\frac{1}{M_r} \bigg)\bigg] \iint_{\q_{\sigma}} u^{r-1}\, dxd\tau \\ \nonumber\nonumber
& \leq &  \gamma \int_{\k_{\sigma \rho}} u^r(x,0)\, dx \\ \nonumber
& & + \frac{\gamma}{(\sigma \rho)^p} \left\{ \sum_{i} \iint_{\q_{\sigma}} u^{r+p_i-2} \, dxd\tau + \sum_i M_r^{p_i-1} \iint_{\q_{\sigma}} u^{r-1} \, dx d\tau \right\}
\end{eqnarray}  
Observe that \eqref{Mr-standard} is a natural assumption: if it is violated then, for some $i\in \{1, \dots, N\}$, then  
\[\bigg( \sup_{0\leq \tau \leq t}\dashint_{\k_{\rho}} u^r(x, \tau) \, dx \bigg)^{\frac{1}{r}} \leq C\rho^{\frac{p}{p_i}} \leq \sum_{k} \bigg( \frac{t}{\rho^p} \bigg)^{\frac{1}{2-p_k}}= \nu_{\Sigma}\] \[\quad \iff \quad \sup_{0\leq \tau \leq t} \int_{\k_{\rho}} u^r(x,\tau)\, dx \leq \sum_k \bigg( \frac{t^r}{\rho^{\lambda_{k,r}}} \bigg)^{\frac{1}{2-p_k}}\]
and \eqref{estimatesLr} is found. Then, as in Theorem \ref{Y}, we estimate the various terms as follows

\begin{equation*}
    \begin{aligned}
        \sum_i \frac{1}{\rho^p} \iint_{\q_{\sigma}} u^{r+p_i-2} \, dxd\tau \leq & \sum_i \bigg( \frac{t}{\rho^p} \bigg) \bigg( \sup_{0\leq \tau \leq t} \int_{\k_{(1+\sigma)\rho}} u^{r+p_i-2}(x,\tau)\, dx \bigg) \\
        & \leq \sum_i \bigg( \frac{t}{\rho^p} \bigg) \bigg( \sup_{0\leq \tau \leq t} \int_{\k_{(1+\sigma)\rho}} u^{r+p_i-2}(x,\tau) \, dx\bigg)^{\frac{p_i+r-2}{r}} (2\rho)^{\frac{N(2-p_i)}{r}}\\
        &= \sum_i \bigg( \frac{t^r}{\rho^{\lambda_{i,r}}} \bigg)^{\frac{1}{r}} \bigg( \sup_{0\leq \tau \leq t} \int_{\k_{(1+\sigma)\rho}} u^r(x,\tau)\, dx \bigg)^{\frac{p_i+r-2}{r}},
    \end{aligned}
\end{equation*}\noindent for $\lambda_{i,r}=N(p_i-2)+pr$, while the second term in the parenthesis of \eqref{estsigma} is managed as follows
\begin{equation*}
    \begin{aligned}
        \sum_i \frac{M_r^{p_i-1}}{\rho^p} \iint_{\q_{\sigma}} u^{r-1} \, dxd\tau \, \,  \leq & \sum_i \bigg( \frac{t}{\rho^p} \bigg)  M_r^{p_i-1}\bigg( \sup_{0\leq \tau \leq t} \int_{\k_{(1+\sigma)\rho}} u^{r}(x,\tau)\, dx \bigg)^{\frac{r-1}{r}} (2\rho)^{\frac{N}{r}} \\
        & \leq \sum_i \bigg( \frac{t^r}{\rho^{\lambda_{i,r}}} \bigg)^{\frac{1}{r}} \bigg( \sup_{0\leq \tau \leq t} \int_{\k_{(1+\sigma)\rho}} u^{r} (x,\tau)\, dx\bigg)^{\frac{p_i+r-2}{r}}.
    \end{aligned}
\end{equation*}

\noindent Plugging these estimates into \eqref{estsigma} we obtain, and applying Young's inequality in each term of the sum,  we get 
\begin{equation}\label{rrrr}
    \begin{aligned}
    \sup_{0\leq \tau \leq t} \int_{\k} u^r(x,\tau) \, dx &\leq \gamma \int_{\k_{(1+\sigma)\rho}} u^r(x,0)\, dx + \sum_i \frac{\gamma}{\sigma^p} \bigg( \sup_{0\leq \tau \leq t} \int_{\k_{(1+\sigma)\rho}} u^r(x,\tau)\, dx \bigg)^{\frac{p_i+r-2}{r}} \bigg(\frac{t^r}{\lambda_{i,r}} \bigg)^{\frac{1}{r}} \\
    &\leq \gamma \int_{\k_{(1+\sigma)\rho}} u^r(x,0)\, dx + \epsilon  \sup_{0\leq \tau \leq t} \int_{\k_{(1+\sigma)\rho}} u^r (x,\tau)\, dx  + \gamma(\epsilon) \sum_i \bigg( \frac{t^r}{\rho^{\lambda_{i,r}}} \bigg)^{\frac{1}{2-p_i}}
    \end{aligned}
\end{equation} \noindent 
From this point on, we perform a standard iteration on $\sigma$: for fixed $\rho>0$ and $n \in \N\cup \{0\}$, we consider the increasing sequence of radii
\[\rho_n:= \rho \sum_{j=0}^n 2^{-j}\ge \rho  \quad \mathrm{so} \ \mathrm{that} \quad \rho_{n+1}=(1+\sigma_n) \rho_n, \ \ \ \mathrm{for} \ \ \sigma_n= \dfrac{\rho_{n+1}-\rho_n}{\rho_n} \geq \frac{1}{2^{n+2}}  \ , \]
by setting \[S_n = \sup_{0\leq \tau \leq t} \int_{\k_{\rho_n}} u^r(x,\tau)\, dx \ ,\]  estimate \eqref{rrrr} now reads
\[S_n \leq \gamma \left\{\epsilon \  S_{n+1}  + b^n \gamma(\epsilon) \left(\int_{\k_{2\rho}} u^r(x,0)\, dx +\sum_i \left(\frac{t^r}{\rho^{\lambda_{i,r}}} \right)^{\frac{1}{2-p_i}}  \right) \right\} \ , \qquad b >1, \]
and the proof is completed once we choose $\epsilon= 1/2b <1$ and let $n \rightarrow \infty$.
\end{proof}

\subsection{Proof of Theorem \ref{elr-decay} concluded}
\begin{proof}
We use \eqref{estimatesLr} to estimate the integral term at the right-hand side of \eqref{estLrlinfty}
\begin{equation*}
    \begin{aligned}
        \sup_{\k_{\rho/2} \times [t/2,t]} u &\leq \gamma t^{-\frac{N}{\lambda_r}} \bigg( \sup_{0\leq \tau \leq t} \int_{\k_{\rho}} \,  u^r(x, \tau)\, dx \bigg)^{\frac{p}{\lambda_r}} + \sum_i \bigg(  \frac{t}{\rho^p}\bigg)^{\frac{1}{2-p_i}}\\
        & \leq \gamma t^{-\frac{N}{\lambda_r}} \bigg(\int_{\k_{\rho}} \, u^r(x,0)\, dx+ \sum_k \bigg(\frac{t^r}{\rho^{\lambda_{k,r}}} \bigg)^{\frac{1}{2-p_k}} \bigg)^{\frac{p}{\lambda_r}} + \sum_i \bigg(  \frac{t}{\rho^p}\bigg)^{\frac{1}{2-p_i}}\\
        & \leq \gamma  t^{-\frac{N}{\lambda_r}} \bigg(\int_{\k_{\rho}} \, u^r(x,0)\, dx\bigg)^{\frac{p}{\lambda_r}}+ \gamma \sum_i \bigg[\bigg(\frac{t}{\rho^p} \bigg)^{\frac{1}{2-p_i}}\bigg]^{\frac{\lambda_{i,r}}{\lambda_r}} + \sum_i \bigg(\frac{t}{\rho^p} \bigg)^{\frac{1}{2-p_i}} \\
    \end{aligned}
\end{equation*}

\end{proof}

\section{Proof of the $L^1$-$L^{\infty}$ estimates}
\label{section-l1linfty}

\subsection*{Intrinsic Geometry. Proof of Theorem \ref{l1-linfty-THM}} 

\begin{proof} We start by considering inequality \eqref{sup-estimate} and then estimate the integral on its right-hand side by \eqref{l1-l1} to get
    \begin{equation*}
        \begin{aligned}
            \|u\|_{\infty, \K_{{\rho}/{2}}(t)\times [t/2,t]} &\leq \gamma t^{-\frac{N}{\lambda}} \bigg(   \inf_{0\leq \tau \leq t}\int_{\K_{2\rho}(t)} u(x,\tau) \, dx + \gamma \bigg( \frac{t}{\rho^{\lambda}} \bigg)^{\frac{1}{2-p}}\bigg)^{\frac{p}{\lambda}}+\gamma \bigg( \frac{t}{\rho^p} \bigg)^{\frac{1}{2-p}} \\&
            \leq \gamma  t^{-\frac{N}{\lambda}} \bigg(   \inf_{0\leq \tau \leq t}\int_{\K_{2\rho}(t)} u(x,\tau) \, dx\bigg)^{\frac{p}{\lambda}}+ \gamma  \bigg( \frac{t}{\rho^p} \bigg)^{\frac{1}{2-p}} .
        \end{aligned}
    \end{equation*}
\end{proof}

\subsection*{Standard Geometry. Proof of Theorem \ref{el1-linfty-THM}}
\begin{proof}

We combine Theorem \ref{lrlinfestimate} with $r=1$ and Theorem \ref{el1l1} to get 
\begin{equation*}
    \sup_{\k_{\rho/2} \times [t/2, t] } u  \leq  \gamma t^{-\frac{N}{\lambda}} \left( \inf_{0\leq \tau \leq t} \int_{2\k_{\rho}} u(x,\tau)  \, dx\right)^{\frac{p}{\lambda}} 
     + \gamma  \sum_i \bigg[\left(\frac{t}{\rho^p}\right)^{\frac{1}{2-p_i}}\bigg]^{\frac{\lambda_i}{\lambda}} + \gamma \sum_i \left(\frac{t}{\rho^p}\right)^{\frac{1}{2-p_i}}.
\end{equation*} 
\end{proof}

\section{Appendix} \label{appendix}

\subsection*{Energy Estimates}
\noindent To the aim of computation, it would be technically convenient to pass from the formulation \eqref{EQ-weak} of local weak solution to its Steklov averaged version, which allows us to perform computations under the integral sign with the approximating functions 
\begin{equation}\label{Steklov}
u_h(x,t)= \begin{cases}
{ \displaystyle \dashint_t^{t+h}u (\cdot, \tau) \, d \tau } , & 0<t <T-h,\quad \text{for} \quad 0<h<T,\\[.8em]
0, & t>T-h,
\end{cases}
\end{equation} defined for all $0<t<T$. This is the same definition as the one presented in  \cite{DB} (see in particular Chapter II for more details), and we refrain from specifying further this procedure, leaving space to what is really new. 

\subsubsection*{Separate Variables Test Functions} \noindent 
\noindent For a compact set $K \subset \Omega$, we will usually test the equation \eqref{EQ-weak}-\eqref{structure-conditions} with functions $\zeta(x) \in C_o^1(K)$ such that
\begin{equation}
    \label{zeta}
    \zeta(x)= \prod_i \zeta_i(x_i)^{p_i},\quad \hat{\zeta}^{j}:= \prod_{i\ne j} \zeta_i(x_i)^{p_i},\qquad 0 \leq \zeta \leq 1,
\end{equation} \noindent with $\zeta_i \in C_o^1(\pi_i(K))$, being $\pi_i$ the euclidean projection to the $i$-th component. Sometimes we will use the notation 
\begin{equation} \label{zetatime} \zeta(x,\tau)=\xi(\tau) \zeta(x),\qquad \qquad \qquad \qquad 0\leq \xi\leq 1,\end{equation}
for $\zeta(x)$ as above and $\xi(\tau)\in C^1_{loc}(0,T)$ a function to be specified at each recurrence. Let $[\tau_1,\tau_2] \subset [0,T]$ be a time interval and $Q=K \times [\tau_1, \tau_2]$ a cylinder inside $\Omega_T$. We denote by 
\[\|\partial_i \zeta \|_{\infty}= \|(\partial_i \zeta) \xi \|_{L^{\infty} (Q)} \quad \text{and} \quad \|\partial_\tau \zeta\|_{\infty}=\|(\partial_\tau \xi)\zeta\|_{L^{\infty}(Q)},  \]
\noindent the essential suprema of $|\partial_i \zeta|$ and $|\partial_{\tau} \zeta|$ in $Q$.

\subsection*{Energy Estimates 1 - Caccioppoli-type Estimates}
\begin{lemma}
Let $u$ be a local weak sub(super)-solution to \eqref{EQ}-\eqref{structure-conditions} and let $ k \in \R$. Let $0\leq \tau_1 <\tau_2\leq T$ and $K \subset \Omega$ be a compact set. Then, there exists a positive constant $\gamma$, depending only on the data, such that for any $\zeta\in C^1_{loc}(0,T; C_{o}^{1}(K))$ of the kind \eqref{zetatime} with $\xi(\tau_1)=0$, we have 
\begin{equation} \label{general-EE}
    \begin{aligned}
        \sup_{\tau_1\leq \tau \leq \tau_2} \int_{K \times \{\tau\}} &(u-k)_+^2 \zeta\, dx+C_o \sum_i \iint_{Q} |\partial_i (u-k)_+ \zeta|^{p_i} dx \, d\tau \\
        &\leq \gamma \sum_i \|\partial_i \zeta_i \|_{\infty}^{p_i} \bigg[1+\bigg(\frac{C}{\|\partial_i \zeta_i \|_{\infty}}\bigg)^{p_i} \bigg]  \iint_{Q} (u-k)_+^{p_i} \,dx d\tau \\
        & + \gamma \|\partial_{\tau} \zeta\|_{\infty } \iint_{Q} \, dxd\tau + \gamma \sum_i C^{p_i} \iint_Q \chi_{[u>k]} \,  dx d\tau,
    \end{aligned}
\end{equation}\noindent where $C\ge 0$ and $C_o > 0$ are the structure constants of \eqref{structure-conditions}.    
\end{lemma}

\begin{proof}
We test equation \eqref{EQ} with $\varphi=(u-k)_+ \zeta$, being $\zeta \in C^1(Q)$ as in \eqref{zetatime}, vanishing on $\partial K$, for all times, and verifying $\zeta(\tau_1,x)=0$, for all $x \in K$. So we  arrive, through a standard Steklov approximation, to 
\begin{equation*}\label{Buddy}
    \begin{aligned}
       \I_1+\I_2:&= \sup_{\tau_1\leq\tau\leq \tau_2}\int_K  \frac{(u-k)_+^2 \zeta}{2} \, dx +\sum_i
 \iint_Q A_i\bigg(\partial_i (u-k)_+ \zeta+ (u-k)_+ (\partial_i \zeta) \bigg) \, dxd\tau  \\
& \qquad \leq  \iint_Q (u-k)_+^2 (\partial_{\tau} \zeta)  \, dx d \tau+ \iint_Q   B (u-k)_+ \zeta \, dxd\tau =: \I_3+ \I_4\, , \end{aligned}
\end{equation*} \noindent being $B,A_i$,for all $i=1,\dots,N$, the Caratheodory functions of \eqref{EQ}-\eqref{structure-conditions}. 
\noindent 
We evaluate the terms separately, using the structure conditions \eqref{structure-conditions} and  Young's inequality \eqref{Young} on each $i$-th term with $q=p_i$, $q'=p_i/(p_i-1)$ to get
\begin{equation*}
    \begin{aligned}
        \I_2& \ge  \sum_i \iint_Q  \bigg(C_o |\partial_i (u-k)_+|^{p_i} -C^{p_i}\chi_{[u>k]} \bigg)\zeta - \bigg(C_1 |\partial_i u|^{p_i-1}+C^{p_i-1} \bigg) (u-k)_+ |\partial_i \zeta_i| p_i\hat{\zeta}^i \zeta_i^{p_i-1}\, dxd\tau\\
        & \ge \sum_i \iint_Q   \bigg(C_o- \gamma \tilde{\epsilon}_i C_1\bigg) |\partial_i (u-k)_+|^{p_i}\zeta- \gamma[\tilde{\gamma}(\tilde{\epsilon}_i) C_1+ 1] (u-k)_+^{p_i} |\partial_i \zeta_i|^{p_i} -  \gamma C^{p_i} \zeta \chi_{[u>k]}\, dxd\tau, 
    \end{aligned}
\end{equation*} \noindent where in the last inequality we have collected the terms
\[|\partial_i \zeta_i|^{p_i} \hat{\zeta}^i= |\partial_i \zeta^{\frac{1}{p_i}}|^{p_i}\leq |\partial_i \zeta_i|, \qquad \text{and} \qquad \hat{\zeta}^i \zeta_i^{p_i}=\zeta\, ,\]
\noindent in order to adjust the powers of $\zeta$.  Again we use Young's inequality for each $i=1, \dots, N$ to estimate
 \begin{equation*}
    \begin{aligned}
        |\I_4|\leq &  \sum_i  \iint_Q   C \bigg(|\partial_i u |^{p_i-1}+ C^{p_i-1} \bigg) (u-k)_+ \zeta \, dxd\tau \\
        & \leq \gamma \sum_i \iint_Q \zeta \epsilon_i |\partial_i (u-k)_+  |^{p_i}  + C^{p_i} (\gamma(\epsilon_i)+1) (u-k)_+^{p_i}  + C^{p_i} \chi_{[u>k]} \, dxd \tau \ .
    \end{aligned}
\end{equation*} \noindent
Choosing suitably $\tilde{\epsilon}_i$ and $\epsilon_i$ small enough for all $i=1, \dots, N$ and joining all the previous estimates together implies, for all $k \in \R$,
\begin{equation*}
    \begin{aligned}
    \sup_{\tau_1 < \tau < \tau_2} & \int_K (u-k)_+^2 \, dx +C_o \sum_i  \iint_Q   \bigg(|\partial_i [(u-k)_+\zeta] |^{p_i} - \gamma (u-k)_+^{p_i}|\partial_i \zeta_i |^{p_i} \bigg) \, dxd\tau\\
& \leq \sup_{\tau_1 < \tau < \tau_2} \int_K (u-k)_+^2 \, dx + \sum_i C_o \iint_Q  |\partial_i (u-k)_+ |^{p_i} \zeta \, dxd\tau\\
&\leq \gamma  \|\partial_{\tau} \zeta \|_{\infty}  \iint_Q (u-k)_+^2 \, dxd\tau + \gamma \sum_i \|\partial_i \zeta_i \|_{\infty}^{p_i} \iint_Q  (u-k)_+^{p_i}  \, dxd\tau\\
& \qquad + \gamma \sum_i \iint_Q  C^{p_i} (u-k)_+^{p_i} \, dxd\tau + \gamma \sum_i \iint_Q  C^{p_i} \chi_{[u>k]} \, dxd\tau.
    \end{aligned}
\end{equation*} \noindent 

\end{proof}

\subsection*{Energy Estimates 2 - Testing with positive powers.}
\begin{lemma} Let $u$ be a non-negative, locally bounded, local weak solution to \eqref{EQ}-\eqref{structure-conditions} satisfying $u \in L^r_{loc}(\Omega)$ for some $r>1$. Let $K_1 \subset  K_2 \subset \Omega$ be compact sets and let $\zeta\in C_o^1(K_2)$ be a cut-off function between $K_1$ and $K_2$ as in \eqref{zeta}. Let $t>0$ be any number such  that the inclusion 
\[Q_j=K_j \times [0,t] \subset \Omega_T,\qquad \qquad \qquad \forall j \in \{1,2\},\]
\noindent is preserved. Then, there  exists a positive constant $\gamma$, depending only on the data, such that 
\begin{equation}\label{parabolic-estimates} 
\begin{aligned}
\sup_{0\leq \tau \leq t} \int_{K_1}& u^r(x, \tau) \, dx \leq \gamma \int_{K_2} u^r(x,0)\, dx +\\
&+ \gamma \sum_{i} \|\partial_i \zeta_i \|_{\infty}^{p_i} \bigg(1+\frac{C^{p_i}}{\|\partial_i \zeta_i\|_{\infty}^{p_i}} \bigg)  \iint_{Q_2} u^{r+p_i-2} \, dxd\tau  \\
&+\gamma\sum_{i} \|\partial_i \zeta_i\|_{\infty}^{p_i} \bigg[\frac{C^{p_i-1}}{\|\partial_i \zeta_i\|_{\infty}^{p_i-1}} +\frac{C^{p_i}}{\|\partial_i \zeta_i \|_{\infty}^{p_i}}\bigg(1+\frac{1}{M_r} \bigg)\bigg] \iint_{Q_2} u^{r-1}\, dxd\tau,
\end{aligned}
\end{equation} being 
\begin{equation}\label{Mr}
   M_r= \bigg( \sup_{0\leq \tau \leq t} \dashint_{K_1} u^r \, dx\bigg)^{1/r}.
\end{equation}
\end{lemma}
\begin{proof}
In the weak formulation \eqref{EQ-weak} choose as a test function, defined over $Q_2$,  \[\varphi=f(u) \zeta= u^{r-1} \bigg(\frac{(u-k)_+}{u}\bigg)^{q}\zeta, \qquad \text{for} \qquad \max\{1, r-1\}<q<r,\] being $\zeta$ as in \eqref{zeta} and $k \in \R^+$ to be determined. We observe that $f(u)=0$ outside the set 
\[[u>k]:=\{ (x,\tau) \in Q_2: \, \, u(x,\tau)>k\}.\]
\noindent Now we define $\displaystyle{F(u)=\int_k^u f(s) \, ds}$ an integral function of $f$ and we observe that 
\begin{equation} \label{trick2}
  (r-1) u^{r-2} \bigg( \frac{(u-k)_+}{u}\bigg)^{q}  \leq f'(u) \leq q u^{r-2} \bigg( \frac{(u-k)_+}{u}\bigg)^{q-1}.
\end{equation}

\noindent The test function $\varphi$ is an admissible one, modulo a Steklov approximation, thanks to the local boundedness of $u$: observe that 
\[\partial_i \varphi = f(u) \partial_i \zeta+ f'(u) \partial_i u \ \zeta \leq \left\{\|\partial_i \zeta\|_{\infty} u^{r-1}+ q \frac{u^{r-1}}{k} |\partial_i u| \, \right\} \chi_{[u>k]}\, \in L^{p_i}_{loc}(\Omega_T).\]
\noindent Passing to the limit the in Steklov approximation, we obtain
\begin{equation*}
    \begin{aligned}
        0= \iint_{\tilde{Q}_2}&\partial_{\tau} F(u)\, \zeta \, dxd\tau + \sum_i  \iint_{\tilde{Q}_2} A_i (\partial_i u) f'(u) \zeta \, dxd\tau 
        \\
        &+ \sum_i \iint_{\tilde{Q}_2} \hat{\zeta}^i f(u) A_i (\partial_i \zeta) \, dxd\tau - \iint_{\tilde{Q}_2} B f(u) \zeta \, dxd\tau=: T_1+T_2+T_3+T_4,
    \end{aligned}
\end{equation*} 
where $\tilde{Q}_2= K_2 \times [0,s]$, for arbitrary $s \in (0,t]$.

\noindent The bound \eqref{trick2} and the fact that $\zeta$ is independent of time allows us to estimate
\[T_1= \int_{K_2} F(u(x,s)) \zeta(x) \, dx-\int_{K_2} F(u(x,0)) \zeta(x) \, dx,\]
\noindent while the structure conditions \eqref{structure-conditions} imply
\begin{equation*}
    \begin{aligned}
        T_2=& \sum_i \iint_{\tilde{Q}_2}  A_i (\partial_i u) f'(u) \chi_{[u>k]}\zeta \, dxd\tau \\
        &\ge \sum_i \iint_{\tilde{Q}_2} \bigg(C_o |\partial_i u|^{p_i} f'(u)-C^{p_i} f'(u)  \bigg)\zeta\, dxd\tau\\
        & \ge \sum_i \iint_{\tilde{Q}_2} \bigg( (r-1)C_o |\partial_i u |^{p_i} u^{r-2} \bigg(\frac{(u-k)_+}{u} \bigg)^{q} -  qC^{p_i} u^{r-2}\bigg(\frac{(u-k)_+}{u} \bigg)^{q-1}     \bigg) \zeta \, dxd\tau\\
        & \ge \sum_i \iint_{\tilde{Q}_2} \bigg( (r-1)C_o |\partial_i u |^{p_i} \frac{f(u)}{u} -  qC^{p_i} u^{r-2}\bigg(\frac{(u-k)_+}{u} \bigg)^{q-1}\bigg) \zeta \, dxd\tau,
    \end{aligned}
\end{equation*}

\[|T_3|\leq \gamma \sum_i\iint_{\tilde{Q}_2} f(u)  \bigg( C_1 |\partial_i u |^{p_i-1} |\partial_i \zeta_i| +C^{p_i-1} |\partial_i \zeta_i|  \bigg) p_i\zeta_i^{p_i-1}\hat{\zeta}^i\, dxd\tau,\]
\[|T_4| \leq \sum_i \iint_{\tilde{Q}_2} \bigg(C |\partial_i u |^{p_i-1} f(u) + C^{p_i} f(u)   \bigg) \zeta \, dxd\tau.\]

\noindent Combining all the estimates we obtain, for all $s \in (0,t]$
\begin{equation*}
    \begin{aligned}
        \int_{K_2} F(u(x,s)) \zeta(x)\, & dx+ (r-1)C_o \sum_i \iint_{\tilde{Q}_2} \frac{f(u)}{u} |\partial_i u |^{p_i} \zeta  \, dxd\tau \leq \int_{K_2} F(u(x,0)) \zeta(x) \, dx\\
        & +\gamma \sum_i \iint_{\tilde{Q}_2} \bigg( C_1+\frac{C}{ |\partial_i \zeta_i|} \bigg) f(u) |\partial_i u |^{p_i-1} |\partial_i \zeta_i |\, \hat{\zeta}^i\zeta_i^{p_i-1}\, dxd\tau \\
        &+ \gamma\sum_i (C^{p_i-1}  \|\partial_i \zeta_i\|_{\infty}) \bigg[1+ \frac{C}{ \|\partial_i \zeta_i\|_{\infty}}  \bigg] \iint_{\tilde{Q}_2} f(u) \, dxd\tau\\
        & + \gamma \sum_i C^{p_i} \iint_{\tilde{Q}_2} u^{r-2} \bigg(\frac{(u-k)_+}{u} \bigg)^{q-1} \zeta \, dxd\tau =: I_1+I_2+I_3+I_4\,.
    \end{aligned}
\end{equation*} \noindent Here we observe that, on the set $[u>k]$, the following holds true 
\[\frac{f(u)}{u} = u^{r-2} \bigg( \frac{(u-k)_+}{u} \bigg)^{q}  \leq  \bigg( \frac{u^{r-1}}{k}\bigg) \quad \qquad \text{and} \qquad f(u) \leq u^{r-1}\quad ,\]
so that we estimate for each $i=1, \dots,N,$ 
\begin{equation*}
    \begin{aligned}
        I_{2,1}= &C_1 \sum_i \iint_{\tilde{Q}_2} f(u) |\partial_i u |^{p_i-1} |\partial_i \zeta_i |\,   \hat{\zeta}^i\zeta_i^{p_i-1}\, dxd\tau \\
        & \leq C_1 \sum_i \epsilon_i \iint_{\tilde{Q}_2} \frac{f(u)}{u} |\partial_i u|^{p_i}  \zeta \, dxd\tau + C_1 \sum_i \gamma(\epsilon_i) \|\partial_i \zeta_i \|_{\infty}^{p_i}\iint_{\tilde{Q}_2} f(u) u^{p_i-1} \zeta \, dxd\tau \\
        & \leq C_1\sum_i \epsilon_i \iint_{\tilde{Q}_2} \frac{f(u)}{u} |\partial_i u|^{p_i} \zeta \, dxd\tau + C_1 \sum_i \gamma(\epsilon_i)  \|\partial_i \zeta_i \|_{\infty}^{p_i} \iint_{\tilde{Q}_2} u^{r+p_i-2} \, \chi_{[u>k]}\, dxd\tau.
    \end{aligned}
\end{equation*} The other integral term does not involve the derivatives of the cut-off function 

\begin{equation*}
    \begin{aligned}
        I_{2,2}= &\sum_i  C \iint_{\tilde{Q}_2}  f(u) |\partial_i u |^{p_i-1} \, \hat{\zeta}^i\zeta_i^{p_i-1}\, dxd\tau \\
        & \leq  \sum_i  \tilde{\epsilon}_i \iint_{\tilde{Q}_2} \frac{f(u)}{u} |\partial_i u|^{p_i} \, \zeta \, dxd\tau + \sum_i C^{p_i}\gamma(\tilde{\epsilon}_i) \iint_{\tilde{Q}_2} f(u) u^{p_i-1} \zeta \, dxd\tau \\
        & \leq  \sum_i  \tilde{\epsilon}_i \iint_{\tilde{Q}_2}\frac{f(u)}{u} |\partial_i u|^{p_i} \zeta \, dxd\tau + \sum_i C^{p_i}\gamma(\tilde{\epsilon}_i)   \iint_{\tilde{Q}_2} u^{r+p_i-2} \, \chi_{[u>k]}\, dxd\tau.
    \end{aligned}
\end{equation*}
\noindent Now we estimate from above $I_3,I_4$ as
\[I_3+I_4 \leq \sum_i  \|\partial_i \zeta_i \|_{\infty}^{p_i} \bigg[ C^{p_i-1}  \|\partial_i \zeta_i \|_{\infty}^{1-p_i}+\frac{C^{p_i}}{ \|\partial_i \zeta_i \|_{\infty}^{p_i}} \bigg( 1+\frac{1}{k} \bigg)\bigg] \iint_{\tilde{Q}_2} u^{r-1} \, \chi_{[u>k]}\, dxd\tau.  \]
Hence, choosing $\epsilon_i$ and $\tilde{\epsilon}_i$ appropriately small, we obtain for all $s \in (0,t]$
\begin{equation}\label{quasi}
    \begin{aligned}
        \int_{K_1} F(u(x,s))\,  dx & \leq \int_{K_2} F(u(x,s)) \zeta(x)\,  dx + \frac{(r-1)C_o}{4} \sum_i \iint_{\tilde{Q}_2} \frac{f(u)}{u}  |\partial_i u |^{p_i} \zeta  \, dxd\tau  \\
        & \leq \int_{K_2} u^r(x, 0)\, dx \\
        & +  \gamma\sum_{i} \|\partial_i \zeta_i\|_{\infty}^{p_i}  \bigg(1+\frac{C^{p_i}}{\|\partial_i \zeta_i \|_{\infty}^{p_i}} \bigg)\iint_{\tilde{Q}_2} u^{r+p_i-2} \, \chi_{[u>k]}\, dxd\tau  \\
        & + \gamma\sum_{i} \|\partial_i \zeta_i\|_{\infty}^{p_i} \bigg[\frac{C^{p_i-1}}{\|\partial_i \zeta_i \|_{\infty}^{p_i-1}} + \frac{C^{p_i}}{\|\partial_i \zeta_i\|_{\infty}^{p_i}}\bigg(1+\frac{1}{k} \bigg)\bigg] \iint_{\tilde{Q}_2} u^{r-1} \, \chi_{[u>k]}\, dxd\tau .
    \end{aligned}
\end{equation} 
since 
\[\int_{K_2} F(u(x, 0))\, dx \leq \int_{K_2} \bigg(\int_0^{u(x, 0)} s^{r-1}ds\bigg)\, dx\, \leq \int_{K_2} u^r(x, 0)\, dx.\]
By choosing $k$ appropriately depending on $M_r$, so that (see for instance \cite{FHV} Prop. 5.1)
\[\sup_{0 \leq \tau \leq t} \dashint_{K_1} u^r (x, \tau)\, dx \leq 2r\bigg( \sup_{0 \leq \tau\leq t} \dashint_{K_1} F(u(x, \tau)) \, dx+ (1+\gamma) k^r |K_1|\bigg) \leq \gamma  \sup_{0 \leq \tau\leq t} \dashint_{K_1} F(u(x, \tau)) \, dx,\] 
estimate \eqref{parabolic-estimates} follows by estimating \eqref{quasi} from below means of this last consideration.

\end{proof}

\begin{remark} The constant $\gamma$ determined along the proof deteriorates as $r \downarrow 1$.    
\end{remark}

\vskip0.2cm \noindent

\subsection*{Energy Estimates 3 - Testing with negative powers}
\begin{lemma}
Let $u$ be a non-negative, local weak super-solution to \eqref{EQ}-\eqref{structure-conditions}. Let $K \subset \Omega$ be a compact set and $0<t<T$ such that $Q= K \times [0,t] \subset \Omega_T$. Then, for all number $\nu>0$ and for all indexes $i=1, \dots, N$ we have the following inequality  
\begin{equation} \label{Negative-EE}
    \begin{aligned}
        \iint_{Q} & \bigg( \sum_j |\partial_j u |^{p_j}\bigg) \tau^{\frac{1}{p_i}} (u+\nu)^{-\frac{2}{p_i}} \zeta\, dxd\tau \leq \gamma t^{\frac{1}{p_i}} \int_{K} (u+\nu)^{\frac{2(p_i-1)}{p_i}}\, dx\\
        & + \gamma  \sum_j  \|\partial_j \zeta_j\|^{p_j}_{\infty}\bigg[1+\bigg( \frac{C}{\|\partial_j \zeta\|_{\infty}}\bigg)^{p_j}
        \bigg]\iint_{Q}  (u+\nu)^{p_j-\frac{2}{p_i}} \tau^{\frac{1}{p_i}}\, dxd\tau \\
        &+ \gamma  \bigg( \sum_j C^{p_j}\bigg)  \iint_{Q} (u+\nu)^{-\frac{2}{p_i}} \tau^{\frac{1}{p_i}}\, dxd\tau,
    \end{aligned}
\end{equation} \noindent for all $\zeta \in C^1(0,t; C^1_o(K))$ of the form \eqref{zetatime}.

\noindent
\end{lemma}

\begin{proof}
   \noindent We test equation \eqref{EQ} repeatedly for $i=1,\dots,N$ with the following test functions
\begin{equation} \label{varfi} \varphi_i(x,\tau)= -\tau^{\frac{1}{p_i}} (u(x,\tau)+\nu)^{1-\frac{2}{p_i}} \zeta(x), \end{equation}\noindent defined in $Q$; where $\zeta$ is a smooth function defined in $K$ of the form \eqref{zeta}. We observe that $\varphi_i(x,0)=0$, for all $x \in K$, and that the function $\varphi_i$, adequately averaged in time, is admissible due to the choice of $\zeta$ and 
\[|\partial_i \varphi_i|\leq \bigg( \frac{2-p_i}{p_i}\bigg) \tau^{\frac{1}{p_i}} \nu^{-2/p_i} \, |\partial_i u| + \tau^{\frac{1}{p_i}} \nu^{(\frac{p_i-2}{p_i})} |\partial_i \zeta| \in L^{p_i}_{loc}(\Omega_T).\]
\noindent In the weak formulation we use Steklov averages (see for instance the monograph \cite{DBGV-mono}) for the interpretation of $\partial_{\tau} u$, to recover by approximation
\begin{equation*}
    \begin{aligned}
        0\ge \int_{K} u \varphi_i \, dx \bigg|_0^t-& \int_0^t \int_{K} u \partial_{\tau} \varphi_i \, dxd\tau + \sum_j \iint_{Q}  A_j  \partial_j \varphi_i \, dxd\tau - \iint_{Q} B \varphi_i \, dxd\tau = I_1-I_2+I_3-I_4.
    \end{aligned}
\end{equation*} As usual in the literature, the parabolic term is estimated by means of Steklov averages thereby getting
\begin{equation*}
    \begin{aligned}
        I_1-I_2 & =\int_{K} u \varphi_i \, dx \bigg|_0^t- \iint_{Q} u (\partial_{\tau} \varphi_i) \, dxd\tau \\
        & =  -\frac{p_i}{2(p_i-1)} \, t^{\frac{1}{p_i}} \int_{K \times\{t\}} (u+\nu)^{\frac{2(p_i-1)}{p_i}} \zeta \, dx  + \frac{1}{2(p_i-1)} \iint_Q (u+\nu)^{\frac{2(p_i-1)}{p_i}} \tau^{\frac{1}{p_i} -1} \zeta \, dxd\tau \\
        & \geq  -\frac{p_i}{2(p_i-1)} \, t^{\frac{1}{p_i}} \int_{K \times\{t\}} (u+\nu)^{\frac{2(p_i-1)}{p_i}} \zeta \, dx
    \end{aligned}
\end{equation*}
\noindent passing to the limit thanks to the condition $u \in C_{loc}(0,T; L^2(K))$, while all the other terms in the Steklov approximation converge to the relative integrals, thanks to the structure conditions and the bound $\nu^{-\alpha}> (u+\nu)^{-\alpha}$, $\nu, \alpha>0$.

\noindent We estimate $I_3$ and $-I_4$ from below by means of Young's inequality
\begin{equation*}
    \begin{aligned}
        I_3 = &\sum_{j}\iint_{Q}  A_j \bigg[\bigg(\frac{2-p_i}{p_i}\bigg) \tau^{\frac{1}{p_i}} (u+\nu)^{-\frac{2}{p_i}} (\partial_j u)\zeta  - \tau^{\frac{1}{p_i}} (u+\nu)^{1-\frac{2}{p_i}} (\partial_j \zeta)  \bigg]\, dxd\tau  \\
        &  \ge \sum_{j} \iint_{Q}  \bigg[C_o|\partial_j u|^{p_j}-C^{p_j} \bigg]\bigg(\frac{2-p_i}{p_i}\bigg) \tau^{\frac{1}{p_i}} (u+\nu)^{-\frac{2}{p_i}}  \zeta \, dx d\tau \\
        &- \sum_{j} \iint_{Q}  \bigg[C_1 |\partial_j u|^{p_j-1} + C^{p_j-1} \bigg] \tau^{\frac{1}{p_i}} (u+\nu)^{1-\frac{2}{p_i}} p_j |\partial_j \zeta_j|\zeta_j^{p_j-1} \hat{\zeta}^j  \, dxd\tau  \\
        & \ge   \sum_{j}  \iint_{Q} \bigg[\bigg(\frac{2-p_i}{p_i}\bigg)C_o-\gamma \epsilon_j C_1 \bigg]|\partial_j u|^{p_j}  \tau^{\frac{1}{p_i}} (u+\nu)^{-\frac{2}{p_i}} \zeta \, dx d\tau \\
        &- \sum_{j} \iint_{Q}  \gamma(\epsilon_j) C_1 (u+\nu)^{p_j-\frac{2}{p_i}} |\partial_j \zeta_j|^{p_j} \tau^{\frac{1}{p_i}}\, dx d\tau \\
        &-\sum_{j} \iint_{Q}  \bigg[\bigg( \frac{2-p_i}{p_i} \bigg)C^{p_j}+ \gamma C^{p_j} \bigg]    (u+\nu)^{-\frac{2}{p_i}}\tau^{\frac{1}{p_i}}  \, dxd\tau\\
        &- \sum_{j} \gamma\iint_{Q}  (u+\nu)^{p_j-\frac{2}{p_i}} |\partial_j \zeta_j|^{p_j}\tau^{\frac{1}{p_i}}
\, dxd\tau .
    \end{aligned}
\end{equation*} 

\begin{equation*}
    \begin{aligned}
        |I_4| \leq& \iint_{Q} \bigg[ \sum_j C \bigg(|\partial_j u|^{p_j-1}+C^{p_j-1} \bigg) \bigg]  (u+\nu)^{1-\frac{2}{p_i}}\tau^{\frac{1}{p_i}}\zeta \, dxd\tau\\
        &\leq \sum_j \iint_{Q}  \bigg[ \tilde{\epsilon}_j |\partial_j u|^{p_j} (u+\nu)^{-\frac{2}{p_i}}\tau^{\frac{1}{p_i}} \zeta  \, dxd\tau+ \tilde{\gamma}(\tilde{\epsilon}_j) C^{p_j} (u+\nu)^{p_j-\frac{2}{p_i}} \tau^{\frac{1}{p_i}} \,  \bigg] dxd\tau\\
        &+\sum_j \iint_{Q} \bigg[ C^{p_j} (u+\nu)^{p_j-\frac{2}{p_i}} \tau^{\frac{1}{p_i}}\zeta + C^{p_j} (u+\nu)^{-\frac{2}{p_i}}\tau^{\frac{1}{p_i}} \zeta \, \bigg] dxd\tau \\& \leq \sum_j \iint_{Q}  \tilde{\epsilon}_j |\partial_j u |^{p_j} (u+\nu)^{-\frac{2}{p_i}} \tau^{\frac{1}{p_i}}\zeta \, dxd\tau \\
        & + \sum_j \iint_{Q}  C^{p_j} \bigg[\tilde{\gamma}(\tilde{\epsilon}_j)+1 \bigg] (u+\nu)^{p_j-\frac{2}{p_i}} \tau^{\frac{1}{p_i}} \, dxd\tau\\
& + \sum_j \iint_Q C^{p_j}(u+\nu)^{-\frac{2}{p_i}} \tau^{\frac{1}{p_i}}  \, dxdt\tau.      
        \end{aligned}
\end{equation*} \noindent 
Now, reabsorbing the terms with $\epsilon_j, \tilde{\epsilon}_j$ on the left-hand side, we obtain 
\begin{equation*} \label{negative-EE}
    \begin{aligned}
        \sum_j \iint_{Q} & |\partial_j u |^{p_j}\tau^{\frac{1}{p_i}} (u+\nu)^{-\frac{2}{p_i}} \zeta\, dxd\tau \leq \gamma t^{\frac{1}{p_i}} \int_{K \times\{t\}} (u+\nu)^{\frac{2(p_i-1)}{p_i}}\zeta\, dx    \\
& + \gamma  \sum_j  \|\partial_j \zeta_j\|^{p_j}_{\infty}\bigg[1+\bigg( \frac{C}{\|\partial_j \zeta_j\||_{\infty}}\bigg)^{p_j}
        \bigg]\iint_{Q}  (u+\nu)^{p_j-\frac{2}{p_i}} \tau^{\frac{1}{p_i}}\, dxd\tau \\
        &+ \gamma \sum_j C^{p_j} \iint_{Q} (u+\nu)^{-\frac{2}{p_i}} \tau^{\frac{1}{p_i}}\, dxd\tau \ .
    \end{aligned}
\end{equation*}

\end{proof}

\begin{remark}
    The constant $\gamma$ deteriorates both as soon as $p_N\uparrow 2$ and as $p_1\downarrow 1$.
\end{remark}

\begin{remark} We observe that all the energy estimates \eqref{general-EE}, \eqref{parabolic-estimates}, \eqref{Negative-EE} recover, when $p_i\equiv p$, known estimates known for the isotropic $p$-Laplacean evolution equations (see for instance the Appendix of \cite{Annali}). This is due to the simple fact that for all $\xi=(\xi_1, \dots, \xi_N) \in \R^N$ there exists an universal constant $\gamma=\gamma(p_i,N)>0$ such that
\[ \frac{1}{\gamma} \sum_i \xi_i^p \leq \| \xi \|^p \leq \gamma \sum_i \xi_i^p, \qquad \text{being} \qquad \|\xi\|= \sqrt{\sum_i \xi_i^2}.  \]
\end{remark}

\vskip0.2cm \noindent 
\subsection*{Algebraic Lemmas} Here we collect two Lemmata evolving sequences of numbers, that can both be found in \cite{DB} (see \cite{CVV} for the anisotropic counterpart), useful along our proofs.

\begin{lemma}{[Fast geometric convergence Lemma]}
\label{fastgeomconv}\vskip0.2cm 

\noindent Let $(Y_n)_n$ be a sequence of positive numbers verifying
\[ Y_{n+1} \leq C b^n \ Y_n^{1+\alpha} \ , \]
being $C>0$, $b>1$ and $\alpha>0$ given numbers. Then the following logical implication holds true
\[Y_o \leq C^{-1/\alpha} \ b^{-1/\alpha^2}\quad \Rightarrow \quad \lim_{n\uparrow \infty} Y_n=0.\]
\end{lemma}

\begin{lemma}{[Iteration Lemma]} \label{iteration} 
\vskip0.2cm 

\noindent If we have a sequence of equibounded numbers $\{Y_n\}$ such that, for constants $\I,b>1$ and $\epsilon \in (0,1)$
\begin{equation}\label{sqn}
    Y_n \leq \epsilon Y_{n+1}+ \I b^n \ , 
\end{equation} 
\noindent then, by a simple iteration, there exists $\gamma >0$ such that 
\[Y_0 \leq \gamma\,  \I.\]
\end{lemma} 

\subsection*{Research Data Policy and Data Availability Statements} All data generated or analysed during this study are included in this article.

\end{document}